\RequirePackage{rotating}
\documentclass[reqno, 12pt]{amsart}
\usepackage[margin=1in]{geometry}
\pdfoutput=1
\usepackage[T2A]{fontenc}
\usepackage[english]{babel}

\usepackage{amssymb,amsthm,amsfonts,amsmath,amscd,bbm,eucal,faktor,xfrac,braket,relsize,nccmath,mathtools,comment,verbatim,wasysym,graphicx,import,float,mathrsfs,color,xparse,chngcntr,mwe,enumitem,tikz,graphicx,wrapfig,lipsum,dsfont,etoolbox,stackrel,bm,upgreek}
\usepackage{url}
\usepackage{genyoungtabtikz}
\usepackage{afterpage}

\usepackage[
backref=true,
backrefstyle=none,
backend=biber,
style=ext-alphabetic,
citestyle=alphabetic,
articlein=false
]{biblatex}

\DefineBibliographyStrings{english}{
  backrefpage = {cited on page},
  backrefpages = {cited on pages},
}

\usepackage[figurename=Figure.]{caption}

\usetikzlibrary{patterns}
\usetikzlibrary{shapes}

\usepackage[
  linktocpage=true
]{hyperref}

\mathtoolsset{showonlyrefs,showmanualtags}

\allowdisplaybreaks

\definecolor{linkblue}{rgb}{0,0.2,0.6}
\definecolor{pictureblue}{rgb}{0,0,1}
\definecolor{pictureblack}{rgb}{0,0,0}
\definecolor{picturered}{rgb}{1,0,0}
\definecolor{picturegreen}{rgb}{0,1,0}

\newcommand{\quot}[2]{{\raisebox{.3em}{$#1$}\bigl/\raisebox{-.3em}{$\!#2$}}}

\newcommand{\RR}{\mathbb{R}_{\geq 0}}

\newcommand\Z{\mathbb Z}
\newcommand\R{\mathbb R}

\newcommand\Pas{\mathbb{P}}

\newcommand\T{\mathcal{T}}

\newcommand\al{\alpha}

\newcommand\Ga{\Gamma}

\newcommand\La{\Lambda}
\newcommand\la{\lambda}

\newcommand\eps{\varepsilon}

\newcommand\wt{\widetilde}
\newcommand\wh{\widehat}

\newcommand\prodd{\prod\limits}
\newcommand\summ{\sum\limits}
\newcommand\limm{\lim\limits}
\newcommand{\ol}{\overline}
\makeatletter
\newcommand*{\isomorphism}{%
  \mathrel{%
    \mathpalette\@isomorphism{}%
  }%
}
\newcommand*{\@isomorphism}[2]{%
  \sbox0{$#1\simeq$}%
  \sbox2{$#1\sim$}%
  \dimen@=\ht0 %
  \advance\dimen@ by -\ht2 %
  \sbox0{%
    \lower2.5\dimen@\hbox{%
      $\m@th#1\relbar\isomorphism@joinrel\rightarrow$%
    }%
  }%
  \rlap{%
    \hbox to \wd0{%
      \hfill\raise\dimen@\hbox{$\m@th#1\sim$}\hfill
    }%
  }%
  \copy0 %
}
\newcommand*{\isomorphism@joinrel}{%
  \mathrel{%
    \mkern-3.4mu %
    \mkern-1mu %
    \nonscript\mkern1mu %
  }%
}
\makeatother

\newcommand\const{\operatorname{const}}

\newcommand\supp{\operatorname{supp}}
\newcommand{\spn}[2]{\operatorname{span}_{#1}\left(#2\right)}
\renewcommand{\Im}{\operatorname{Im}}

\DeclareMathOperator{\K}{K_0}
\DeclareMathOperator{\KK}{K_0^+}

\newcommand\zig{\EuScript{Z}}
\newcommand\qsym{\mathit{QSym}}

\newcommand{\bw}{\mathsf{bw}}
\newcommand{\z}{\mathsf{z}}
\newcommand{\fl}{\mathsf{fl}}

\newcommand{\lr}[1]{\left( #1 \right)}
\newcommand{\biglr}[1]{\Bigl( #1 \Bigr)}

\newcommand{\p}[1]{\overset{\raisebox{-0.5ex}{$\scriptstyle #1$}}{+}}
\newcommand{\m}[1]{\overset{\raisebox{-0.5ex}{$\scriptstyle #1$}}{-}}
\newcommand{\pp}{\overset{\infty}{+}}
\newcommand{\mm}{\overset{\infty}{-}}
\newcommand{\ml}[1]{\mathlarger{#1}}
\newcommand{\ms}[1]{\mathsmaller{#1}}

\renewcommand{\binom}[2]{\begin{pmatrix}
#1\\ 
#2
\end{pmatrix}
}
\newtheorem{theorem}{Theorem}
\newtheorem{proposition}[theorem]{Proposition}

\newtheorem{lemma}[theorem]{Lemma}

\newtheorem*{theorem*}{Theorem}
\newtheorem*{proposition*}{Proposition}
\newtheorem*{corollary*}{Corollary}

\theoremstyle{definition}
\newtheorem{definition}[theorem]{Definition}
\newtheorem{definition-proposition}{Definition-Proposition}
\newtheorem{remark}[theorem]{Remark}
\newtheorem{observation}[theorem]{Observation}
\newtheorem{notation-definition}[theorem]{Notation-Definition}
\newtheorem{example}[theorem]{Example}

\newtheorem*{remark*}{Remark}
\newtheorem*{observation*}{Observation}
\newtheorem*{example*}{Example}

\numberwithin{definition-proposition}{section}
\numberwithin{theorem}{section}

\title{Semifinite harmonic functions on the Zigzag graph}

\author{Nikita Safonkin}

\address{Skolkovo Institute of Science and Technology, Moscow, Russia \& National Research University Higher School of Economics, Moscow, Russia.}

\email{safonkin.nik@gmail.com}

\begin{document}

\begin{abstract}
    We study semifinite harmonic functions on the zigzag graph, which corresponds to Pieri's rule for the fundamental quasisymmetric functions $\{F_{\la}\}$. The main problem, which we solve here, is to classify the indecomposable semifinite harmonic functions on this graph. We describe the set of classification parameters and an explicit construction that produces a semifinite indecomposable harmonic function out of every point of this set. We also establish a semifinite analog of the Vershik-Kerov ring theorem.
\end{abstract}

    \maketitle

    \setcounter{tocdepth}{2}
    \tableofcontents

\section{Introduction}

The zigzag graph $\zig$ is a $\mathbb{Z}_{\geq0}$-graded graph whose vertices correspond to the compositions of natural numbers and which edges are defined via the Pieri rule for fundamental quasisymmetric functions. The zigzag graph is the Hasse diagram of the subword order on the set of binary words, see Section \ref{par zig}.

A function on the set of vertices of the zigzag graph is called \textit{harmonic}, if the value at any vertex $\la$ from the $n$-th level equals the sum of values at all vertices from the $n+1$-st level joined with $\la$ by an edge. Harmonic functions are assumed to be real-valued and non-negative; furthermore, they can take the value $+\infty$. A. V. Gnedin and G. I. Olshanski have described all \textit{finite} harmonic functions on the zigzag graph \cite{gnedin_olsh2006}. These functions are in bijection with the probability measures on the space of pairs of disjoint open subsets of $(0,1)$. Moreover, the indecomposable harmonic functions correspond to the delta-measures.

In the present paper we describe all indecomposable \textit{semifinite} harmonic functions on the zigzag graph. The semifiniteness condition means that at some vertices the function takes the value $+\infty$; for the precise definition, see Section \ref{sec gen theory}. The main result is Theorem \ref{main theorem}. It shows that the parameters of this classification are the so-called \textit{semifinite zigzag growth models}. A semifinite zigzag growth model is a pair of an infinite zigzag, which is a formal collection of possibly infinite rows and columns, and a tuple of real positive numbers that are summed to $1$ and attached to infinite rows and columns. These numbers are treated as growth frequencies of the corresponding rows and columns. See Definitions \ref{framed oriented pb} and \ref{332}, and Remark \ref{rem zero temp}. One can treat a semifinite zigzag growth model as a pair of disjoint open subsets of $(0,1)$ consisting of finite number of intervals plus some additional discrete parameters (a tuple of zigzags assigned to boundary points of the intervals) by viewing each growth frequency as length of a subinterval of $(0,1)$, see Remark \ref{rem456}. Forgetting about the discrete part of the data, we obtain a pair of disjoint open subsets of $(0,1)$. It turns out that this pair corresponds to a finite indecomposable harmonic function that appears in the semifinite analog of the Vershik-Kerov ring theorem, see Theorem \ref{multiplicativity theorem} and Proposition \ref{prop mult for zig}.

With the help of the ergodic method, A. M. Vershik and S. V. Kerov have described the semifinite harmonic functions on the Young and Kingman graphs \cite{kerov_example,vershik_Kerov83}. This method is applicable to any branching graph, but it requires an evaluation of some limit, which turns out to be complicated for the zigzag graph. Instead of the ergodic method, we will follow another approach proposed by A. J. Wassermann \cite[chapter III, Section 6]{wassermann1981}, see also \cite{safonkin21}.

Note that the zigzag graph $\zig$ is the Bratteli diagram of some AF-algebra, and indecomposable semifinite harmonic functions on it correspond to normal factor representations of type I$_\infty$ and II$_\infty$ of this algebra.

\subsection{Organization of the paper}
In Section \ref{par zig} we introduce the zigzag graph. In Section \ref{sec gen theory} the definition of semifinite harmonic functions is given and Wassermann's method is briefly discussed. Section \ref{sec coideals zig} is devoted to saturated primitive coideals of the zigzag graph. We find out which of them do not admit strictly positive semifinite indecomposable harmonic functions. Inside each of the remaining coideals we introduce an ideal, which we use in the subsequent section. Furthermore, we discuss some of the properties of that ideal and give a couple of examples. Main results of the paper, Theorem \ref{main theorem} and Proposition \ref{prop mult for zig}, are discussed in Section \ref{main section}. 

\subsection{Acknowledgements}
I am deeply grateful to Grigori Olshanski for many useful comments and stimulating discussions. I would like to thank Pavel Nikitin for careful reading of the paper and helpful discussions. Supported in part by the Simons Foundation. Partially supported by the Basic Research Program at the HSE University.

\section{The zigzag graph}\label{par zig}
In this section we recall a few notions on the zigzag graph \cite{gnedin_olsh2006}, see also \cite{tarrago}.

Let us consider compositions (ordered partitions) of natural numbers. We identify them with the ribbon diagrams, which are connected skew Young diagrams written in the French notation and containing no $2\times 2$ blocks of boxes. A composition $\la=(\la_1,\ldots,\la_l)$ is identified with the ribbon Young diagram having $\la_i$ boxes in the $i$-th row. For instance, the only one composition of $1$ gets identified with $\square$. The number of boxes in $\la$ equals $|\la|=\la_1+\ldots+\la_l$. We treat ribbon Young diagrams as zigzags crawling from the top-left corner to the bottom-right corner. There is a bijection between the zigzags and the binary words, which we will discuss in details. 

\textit{A binary word} is a word in the alphabet of two symbols, $+$ and $-$. We will use the following conventions $$\p{n}=\underbrace{+\ldots+}_{n}\ \ \ \ \ \text{and}\ \ \ \ \ \m{n}=\underbrace{-\ldots-}_{n}.$$ 

The bijection between the zigzags and the binary words is as follows. From left to right we read the symbols off the binary word and add boxes to the simplest zigzag $\square$. If the symbol is $+$, then we add a box in the horizontal direction to the right, and if the symbol is $-$, then we add a box in the vertical direction to the bottom. For instance, the binary word $-+$ corresponds to the zigzag with one box in the first row and two boxes in the second row. The binary word corresponding to a zigzag $\la$ will be denoted by $\bw(\la)$. So, $\bw(\square)$ is the empty binary word. The number of symbols in a binary word will be denoted by $|\cdot|$, so $|\bw(\la)|=|\la|-1$.

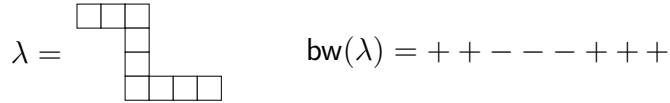
\begin{figure}[h]
    \centering
    \begin{tikzpicture}[yscale=1.0,xscale=1.0]
    \node at (2,8){ \scalebox{0.7}{$\gyoung(;;;,::;,::;,::;;;;)$}};
    \node at (0.5,8){$\la=$};
    \node at (6.5,8){$\bw(\la)=++---++\,+$};
    \end{tikzpicture}
    \caption{A zigzag diagram and its binary word}
\end{figure}

For binary words $a$ and $b$ we write $a\nearrow b$ if and only if $|b|=|a|+1$ and $a$ can be obtained from $b$ by deleting a single symbol. For zigzags $\la$ and $\mu$ we write $\la\nearrow \mu$ if and only if $\bw(\la)\nearrow \bw(\mu)$.

\begin{definition}
The zigzag graph $\zig$ is a graded graph $\zig=\bigsqcup\limits_{n\geq 0} \zig_n$, where $\zig_n$ is the set of all zigzags consisting of $n$ boxes. By definition $\zig_0$ is a singleton $\zig_0=\{\diameter\}$. There is an edge going from $\la$ to $\mu$ if and only if $\la\nearrow \mu$. All edges of $\zig$ are by definition simple.
\end{definition}

\begin{figure}[h]
    \centering
    \scalebox{0.8}{
\begin{tikzpicture}[yscale=1.0,xscale=1.0]

\node (0) at (5,0){$\diameter$};

\node (11) at (5,1){$\gyoung(;)$};

\node (21) at (7.5,2){$\gyoung(;;)$};
\node (22) at (2.5,2){$\gyoung(;,;)$};

\node (31) at (9,3.5){$\gyoung(;;;)$};
\node (32) at (6,3.5){$\gyoung(;;,:;)$};
\node (33) at (4,3.5){$\gyoung(;,;;)$};
\node (34) at (1,3.5){$\gyoung(;,;,;)$};

\node (41) at (12,6.5){$\gyoung(;;;;)$};
\node (42) at (10,6.5){$\gyoung(;;;,::;)$};
\node (43) at (8,6.5){$\gyoung(;;,:;;)$};
\node (44) at (6,6.5){$\gyoung(;;,:;,:;)$};
\node (45) at (4,6.5){$\gyoung(;,;;;)$};
\node (46) at (2,6.5){$\gyoung(;,;;,:;)$};
\node (47) at (0,6.5){$\gyoung(;,;,;;)$};
\node (48) at (-2,6.5){$\gyoung(;,;,;,;)$};

\draw (0) to (11); 
\draw (11) to (21);
\draw (11) to (22);
\draw (21) to (31);
\draw (21) to (32);
\draw (21) to (33);
\draw (22) to (32);
\draw (22) to (33);
\draw (22) to (34); 
\draw (31) to (41);
\draw (31) to (42);
\draw (31) to (43);
\draw (31) to (45);
\draw (32) to (42);
\draw (32) to (43);
\draw (32) to (44);
\draw (32) to (46);
\draw (33) to (43);
\draw (33) to (45);
\draw (33) to (46);
\draw (33) to (47);
\draw (34) to (44);
\draw (34) to (46);
\draw (34) to (47);
\draw (34) to (48);

\end{tikzpicture}
}
    \caption{The first few levels of the zigzag graph $\zig$}
\end{figure}
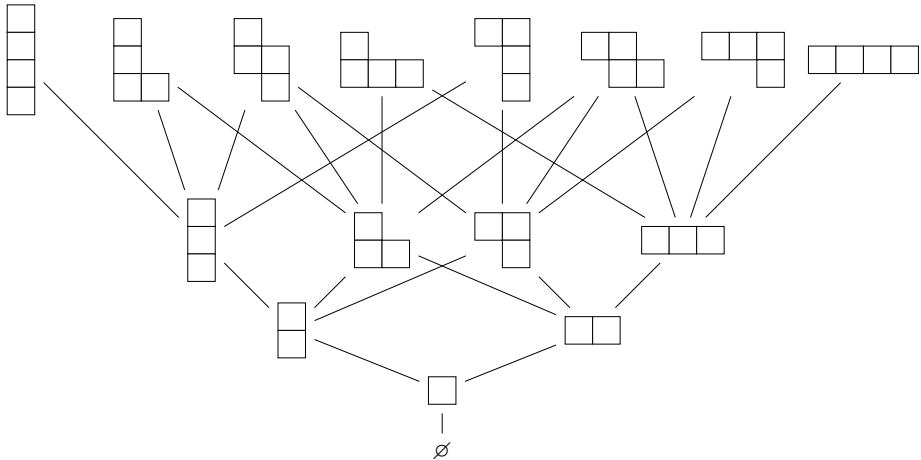

Let $\qsym$ be the algebra of quasisymmetric functions and $\{F_{\la}\}_{\la\in\zig}$ be the fundamental quasisymmetric functions defined in \cite{gessel_84}, see also \cite[p. 357]{stanley2}. There is a Pieri-type rule for this basis $$F_{\ms{\square}}F_{\la}=\summ_{\mu:\la\nearrow \mu}F_{\mu},$$ 
which reflects the branching rule for the zigzag graph, see \cite[p. 482, Exercise 7.93]{stanley2} or \cites[p. 35, (3.13)]{quaisymmetric_book}. 

\section{Branching graphs and semifinite harmonic functions}\label{sec gen theory}

In this section we give the definition of semifinite harmonic functions and briefly recall their main properties. The proof of Theorem \ref{main theorem} relies on Proposition \ref{co-ideal}, Theorem \ref{theorem ext}, Proposition \ref{prop prim coideal}, and Theorem \ref{theorem mult graph}. A similar approach was used by A. Wassermann for the description of semifinite harmonic functions on the Young graph \cite[chapter III, Section 6]{wassermann1981}. The material from the present section is contained \footnote{The papers \cite{strat_voic1975,Bratteli1972,wassermann1981} deal with operator algebras, while we work with combinatorial objects and talk about harmonic functions on branching graphs instead of traces on AF-algebras.} in \cite{strat_voic1975,wassermann1981,Bratteli1972} and \cite{vershik_Kerov83,versh_ker_85,kerov_vershik1990}. A detailed combinatorial exposition can be found in \cite{safonkin21}.

\subsection{Coideals of branching graphs}\label{br beginning}
In all what follows, graded graphs are assumed to have finite levels and simple edges, which join only some vertices from adjacent levels. We assume that the edges are oriented from lower to higher levels and each vertex is joined by an edge with some vertex from the higher level. The level of a vertex $\la$ is denoted by $|\la|$. We write $\la\nearrow \mu$, if there is an edge going from $\la$ to $\mu$. \textit{A path} is a (finite or infinite) sequence of vertices  $\la_1,\la_2,\la_3,\ldots$ such that $\la_i\nearrow \la_{i+1}$ for every $i$. We say that $\nu$ lies above $\mu$ if $\nu$ belongs to a higher level than $\mu$ and they can be connected by a path. In this case, we write $\nu>\mu$. The number of paths going from $\la$ to $\mu$ is denoted by $\dim(\la,\mu)$ and is called the shifted dimension. A graded graph $\Ga$ is called \textit{a branching graph}, if $\Ga_0=\{\diameter\}$ is a singleton and for any vertex there is a vertex from the lower level joined with that by an edge.

\begin{definition}
    A subset of vertices $I$ of a graded graph $\Ga$ is called an \textit{ideal}, if for any vertices $\la\in I$ and $\mu\in\Ga$ such that $\mu>\la$ we have $\mu\in I$. A subset $J\subset \Ga$ is called a \textit{coideal}, if for any vertices $\la\in J$ and $\mu\in\Ga$ such that $\mu<\la$ we have $\mu\in J$.
\end{definition}

There is a bijective correspondence $I\leftrightarrow\Ga\backslash I$ between ideals and coideals. Let $J$ be a coideal and $I=\Ga\backslash J$ be the corresponding ideal. Then the following conditions are equivalent:
	\begin{enumerate}[label=\theenumi)]
	    \item if $\left\{\mu \;\middle|\; \la\nearrow\mu\right\}\subset I$, then $\la\in I$
	    
	    \item for any $\la\in J$ there exists a vertex $\mu\in J$ such that $\la\nearrow \mu$.
	\end{enumerate}
	
	\begin{definition}
	An ideal $I$ and the corresponding coideal $J$ are called \textit{saturated}, if they satisfy the conditions above. A saturated ideal $I$ is called \textit{primitive}, if for any saturated ideals $I_1,I_2$ such that $I=I_1\cap I_2$ we have $I=I_1$ or $I=I_2$. A saturated coideal $J$ is called \textit{primitive}, if for any saturated coideals $J_1,J_2$ such that $J=J_1\cup J_2$ we have $J=J_1$ or $J=J_2$.
	\end{definition}
	
	The bijection $I\leftrightarrow\Ga\backslash I$ maps primitive saturated ideals to primitive saturated coideals and vice versa. 
	
	Let $\Ga$ be a branching graph. The space of infinite paths in $\Ga$ starting at $\diameter$ will be denoted by $\T(\Ga)$. To every path $\tau=\left(\diameter\nearrow\la_1\nearrow \la_2\nearrow \ldots\right)\in\T(\Ga)$ we associate the saturated primitive coideal $\Ga_{\tau}=\bigcup\limits_{n\geq 1}\{\la\in\Ga\mid \la\leq \la_n\}$.
	
	\begin{proposition}\cites{Bratteli1972,strat_voic1975}[p.129]{wassermann1981}[Proposition 2.6]{safonkin21}\label{co-ideal}
	\begin{enumerate}[leftmargin=5ex, label=\theenumi)]
	    \item A saturated coideal $J$ of a graded graph is primitive if and only if for any two vertices $\la_1,\la_2\in J$ we can find a vertex $\mu\in J$ such that $\mu\geq\la_1,\la_2$.
	    
	    \item Every saturated primitive coideal of a branching graph is of the form $J=\Ga_{\tau}$ for some path $\tau\in \T(\Ga)$.
	\end{enumerate}
	\end{proposition}

\begin{definition}
A graded graph $\Ga$ is called \textit{primitive} if for any vertices $\la_1,\la_2\in \Ga$ there exists a vertex $\mu\in \Ga$ such that $\mu\geq\la_1,\la_2$.
\end{definition}

\subsection{Semifinite harmonic functions}\label{qwer}

Let $\Ga$ be a graded graph. 

\begin{definition}
 A function $\varphi\colon \Ga\rightarrow \R_{\geq 0}\cup\{+\infty\}$ is called \emph{harmonic}, if it satisfies the following condition 
 $$\varphi(\la)=\summ_{\mu:\la\nearrow \mu}\varphi(\mu).$$ 
 The set of all vertices $\la\in\Ga$ with $\varphi(\la)<+\infty$ is called the \textit{finiteness ideal} of $\varphi$. We denote the \textit{zero ideal} $\left\{\la\in\Ga\; \middle|\; \varphi(\la)=0\right\}$ by $\ker{\varphi}$ and the \textit{support} $\{\la\in\Ga\mid \varphi(\la)>0\}$ by $\supp{\varphi}$.
\end{definition}

Note that there is an obvious bijection between harmonic functions on a graded graph $\Ga$ with the given support $J$ and strictly positive harmonic functions on $J$. We will use it many times throughout the paper.

The symbol $\K(\Ga)$ stands for the $\R$-vector space spanned by the vertices of $\Ga$ subject to the following relations $$\la=\summ_{\mu:\la\nearrow \mu}\mu,\ \ \forall \la\in\Ga.$$ 
The symbol $\KK(\Ga)$ denotes the positive cone in $\K(\Ga)$, generated by the vertices of $\Ga$, i.e. $\KK(\Ga)=\spn{\RR}{\la\mid\la\in\Ga}$. The partial order, defined by the cone $\KK(\Ga)$, is denoted by $\geq_K$. That is $a\geq_K b\iff a-b\in \KK(\Ga)$. For instance, if $\la\geq \mu$, then $\mu\geq_K\dim(\mu,\la)\cdot\la$.

\begin{remark}
Notation $\K(\Ga)$ is motivated by the fact that the vector space $\K(\Ga)$ can be identified with the Grothendieck $\K$-group of the corresponding AF-algebra. Under such a bijection the cone $\KK(\Ga)$ gets identified with the cone of true modules \cite[Theorem 13 on page 32]{versh_ker_85}. 
\end{remark}

The $\RR$-linear map $\KK(\Ga)\rightarrow \RR\cup\{+\infty\}$, defined by a harmonic function $\varphi$, will be denoted by the same symbol $\varphi$. 

\begin{definition}\label{semifinite def}
A harmonic function $\varphi$ is called \emph{semifinite}, if it is not finite and for any $a\in\KK(\Ga)$ the map $\varphi\colon\KK(\Ga)\rightarrow \RR\cup\{+\infty\}$ enjoys the following property 
$$\varphi(a)=\sup\limits_{\substack{b\in\KK(\Ga)\colon  b\leq_K a,\\ \varphi(b)<+\infty}}\varphi(b).$$
\end{definition}

\begin{remark}\label{remark 321}
 A harmonic function $\varphi$ is semifinite if and only if there exists an element $a\in\KK(\Ga)$ with $\varphi(a)=+\infty$ and for any such $a$ we can find a sequence $\{a_n\}_{n\geq 1}\subset \KK(\Ga)$ such that
 \begin{itemize}
     \item $a_n\leq_K a$,
     
     \item $\varphi(a_n)<+\infty$,
     
     \item $\lim\limits_{n\to+\infty}\varphi(a_n)=+\infty$.
 \end{itemize} 
 We will call this $\{a_n\}_{n\geq 1}$ \textit{an approximating sequence}.
\end{remark}

\begin{definition}\label{def indec}
A semifinite harmonic function $\varphi$ is called \textit{indecomposable}, if for any finite or semifinite harmonic function $\varphi'$ which does not vanish identically on the finiteness ideal of $\varphi$ and satisfies the inequality $\varphi'\leq \varphi$ we have $\varphi'=\const\cdot\varphi$ on the finiteness ideal of $\varphi$.
\end{definition}
 
\begin{theorem}\cites[Theorem 7 стр.143, Corollary стр.144 ]{wassermann1981}[Theorem 3.14]{safonkin21}\label{theorem ext}
Let $I$ be an ideal of a primitive graded graph $\Ga$. The strictly positive indecomposable finite and semifinite harmonic functions on $\Ga$ are in a bijective correspondence with the similar functions on $I$. This bijection is defined by the restriction of functions on $\Ga$ to the ideal $I$.
\end{theorem}

\begin{proposition}\label{prop prim coideal}\cites[p.35 Lemma 12]{versh_ker_85}[Proposition 3.16]{safonkin21}
	Let $\Ga$ be a graded graph and $\varphi$ be an indecomposable finite or semifinite harmonic function on it. Then $\supp(\varphi)$ is a primitive coideal.
\end{proposition}

\subsubsection{Semifinite harmonic functions on multiplicative graphs}
\begin{definition}\label{multiplicative def}
A branching graph $\Ga$ is called \textit{multiplicative}, if there exists an associative $\Z_{\geq 0}$-graded $\R$-algebra $A=\bigoplus\limits_{n\geq0} A_n$, $A_0=\R$ with a distinguished basis of homogeneous elements $\{a_{\la}\}_{\la\in\Ga}$ such that
\begin{enumerate}[label=\theenumi)]
    \item $\deg{a_{\la}}=|\la|$
    
    \item $a_{\diameter}$ is the identity in $A$
    
    \item $\wh{a}\cdot a_{\la}=\summ_{\mu:\la\nearrow\mu} a_{\mu}$ for $\wh{a}=\summ_{\nu\in\Ga_1}a_{\nu}$ and any vertex $\la\in\Ga$. 
\end{enumerate}

Moreover, we assume that the structure constants of $A$ with respect to the basis $\{a_{\la}\}_{\la\in\Ga}$ are non-negative. 
\end{definition}

Recall the Vershik-Kerov ring theorem.

We say that a harmonic function $\varphi$ on a branching graph $\Ga$ is \textit{normalized} if $\varphi(\diameter)=1$.

\begin{theorem}\cites[Theorem p.134]{vershik_Kerov83}[Proposition 8.4]{gnedin_olsh2006}\label{vershik kerov ring theorem}
A finite normalized harmonic function $\varphi$ on a multiplicative branching graph $\Ga$ is indecomposable iff the corresponding functional on $A$ is multiplicative: $\varphi\lr{a\cdot b}=\varphi\lr{a}\cdot \varphi\lr{b}\ \forall a,b\in A$.
\end{theorem}

The following semifinite analog of the ring theorem holds.

\begin{theorem}\cites[Theorem p.144]{vershik_Kerov83}[Theorem 4.4]{safonkin21}\label{multiplicativity theorem}
For any semifinite indecomposable harmonic function $\varphi$ on a multiplicative branching graph $\Ga$ there exists a finite normalized indecomposable harmonic function $\psi$ such that for any $\mu\in \Ga$ with $\varphi(\mu)<+\infty$ we have $\varphi(a_{\la}\cdot a_{\mu})=\psi(a_{\la})\cdot\varphi(a_{\mu})$.
\end{theorem}

\begin{theorem}\cites[Theorem 8 p.146]{wassermann1981}[Theorem 4.5]{safonkin21}\label{theorem mult graph}
Let $\Ga$ be a multiplicative graph. If $a_{\la}a_{\mu}\neq0$ for any $\la,\mu\in\Ga$, then the graph $\Ga$ admits no strictly positive semifinite indecomposable harmonic functions.
\end{theorem}

\begin{proposition}\cites[p. 371, the paragraph just before Theorem 3.5]{boyer_symplectic}[Corollary 4.6]{safonkin21}\label{cunning Boyer}
If a multiplicative branching graph $\Ga$  admits a strictly positive indecomposable finite harmonic function, then $\Ga$ possesses no strictly positive semifinite indecomposable harmonic functions. 
\end{proposition}

\subsubsection{Boyer's lemma}\label{boyer's lemma}
Let $\Ga_1$ and $\Ga_2$ be graded graphs with a given graded graph isomorphism $\Ga_1\rightarrow \Ga_2,\ \la\mapsto \la'$. Let $\Ga$ be another graded graph such that $\lr{\Ga}_0=\lr{\Ga_1}_0$, for $n\geq 1$ $\lr{\Ga}_n=\lr{\Ga_1}_n\sqcup \lr{\Ga_2}_{n-1}$, and the oriented edges $\la\nearrow \mu$ of $\Ga$ are of the following three kinds:
\begin{itemize}
    \item $\la,\mu\in\Ga_1$ and $\la\nearrow\mu$;
    
    \item $\la,\mu\in\Ga_2$ and $\la\nearrow\mu$;
    
    \item $\la\in\Ga_1$ and $\mu=\la'\in\Ga_2$.
\end{itemize}

So, we allow some edges going from $\Ga_1$ to $\Ga_2$, but not vice versa.

\begin{lemma}\cites[Boyer's Lemma p.149]{wassermann1981}[Corollary 5.5, Remark 5.6]{safonkin21}\label{rem boyers classic}
Let $\la\in \Ga_1$ and $\varphi$ be a harmonic function on $\Ga$ such that $\varphi(\la')>0$. Then $\varphi(\la)=+\infty$.
\end{lemma}

\subsubsection{Generalized Boyer's lemma}

\begin{lemma}\cites[Theorem 1.10]{boyer}[Proposition 5.3]{safonkin21}\label{gen boyer's lemma}
Let $\Ga$ be a graded graph and $\varphi$ be a harmonic function on it. Assume that $I\subset \Ga$ is an ideal, $J=\Ga\backslash I$ is the corresponding coideal and we are given a fixed vertex $\la\in J_n$ lying on the $n$-th level of $J$. Suppose that there exists a vertex $\la'$ such that $\varphi(\la')>0$ and for any large enough $l$ and any vertex $\eta\in I_{n+l+1}$ the following inequality holds
\begin{equation}\label{ineq gen boyers}
\summ_{\mu\in J_{n+l},\mu\nearrow\eta}\dim(\la,\mu)\geq \dim(\la',\eta).
\end{equation}

Then $\varphi(\la)=+\infty$.  
\end{lemma}

\begin{remark}
Lemma \ref{gen boyer's lemma} does not look exactly like \cite[Proposition 5.3]{safonkin21}, but in fact they are almost equivalent.
\end{remark}

\subsubsection{Direct product of branching graphs}\label{sect prod graph}

\begin{definition}\label{direct product graph def}
By the \textit{direct product} of branching graphs $\Ga_1$ and $\Ga_2$ we mean the branching graph $\Ga_1\times \Ga_2$, where
$$\lr{\Ga_1\times\Ga_2}_k=\bigsqcup\limits_{\substack{n,m\geq0: \\ n+m=k}}\lr{\Ga_1}_n\times\lr{\Ga_2}_m$$
and $(\la_1,\mu_1)\nearrow(\la_2, \mu_2)$ if and only if one of the following holds
\begin{itemize}
    \item $\la_1=\la_2$, $\mu_1\nearrow\mu_2$;
    
    \item $\la_1\nearrow\la_2$, $\mu_1=\mu_2$.
\end{itemize}
\end{definition}

\begin{proposition}\cite[Proposition A.4]{safonkin21}\label{product graph}
Let $\Ga_1,\ldots, \Ga_n$ be some branching graphs and $\varphi$ be a finite strictly positive indecomposable normalized harmonic function on the graph $\Ga_1\times\ldots\times\Ga_n$. Then there exist finite normalized strictly positive indecomposable harmonic functions $\varphi_1,\ldots, \varphi_n$ on graphs $\Ga_1,\ldots,\Ga_n$ and a tuple of positive real numbers $w_1,\ldots,w_n$ with $w_1+\ldots+w_n=1$ such that for any $\la_1\in\Ga_1,\ \ldots\ ,$ $\la_n\in\Ga_n$ we have
$$\varphi(\la_1,\ldots,\la_n)=w_1^{|\la_1|}\ldots w_n^{|\la_n|}\varphi_1(\la_1)\ldots \varphi_n(\la_n).$$ 
Moreover, these $\varphi_1,\ldots,\varphi_n$ and $w_1,\ldots,w_n$ are uniquely defined. 
\vspace{1\baselineskip}

\end{proposition}

\section{Coideals of the zigzag graph}\label{sec coideals zig}

The algebra $\qsym$ contains no zero divisors, since it is a subalgebra of the formal power series algebra in countably many variables. Then Theorem \ref{theorem mult graph} implies that $\zig$ possesses no strictly positive indecomposable semifinite harmonic functions. From Proposition \ref{prop prim coideal} it follows that the support of any indecomposable semifinite harmonic function on $\zig$ is a primitive saturated coideal. In this section we explicitly describe all primitive saturated coideals of $\zig$. Furthermore, we specify the coideals corresponding to the supports of the finite indecomposable harmonic functions. By Proposition \ref{cunning Boyer} none of these coideals can be realised as the support of an indecomposable semifinite harmonic function.

\subsection{Saturated coideals of the zigzag graph}

Each binary word can be uniquely represented as a consecutive union of \textit{blocks} with alternating signs. By a \textit{block} we mean a tuple of symbols of the same sign. For instance, the word $+-\p{3}$ splits into three blocks, $+$, $-$, and $\p{3}$. So, a block can be positive or negative depending on the sign of symbols. As for zigzags, these positive and negative blocks correspond to rows and columns.

\begin{definition}\label{332}
By a \textit{cluster} we mean a symbol, $+$ or $-$, with an assigned to it formal positive multiplicity, which may be infinite. We say that a cluster is \textit{infinite}, if its multiplicity is infinite, otherwise we say that the cluster is \textit{finite}. A \emph{template} is an ordered collection of alternating clusters. Furthermore, we always assume that a template contains at least one infinite cluster. 
\end{definition}

For instance, $\pp\mm\p{3}\mm$ is a template while $\p{3}\m{2}$ is not.

Each template can be thought of as an infinite zigzag consisting of finite number of possibly infinite rows and columns. Infinite rows and columns correspond to infinite clusters of this template. The infinite zigzag corresponding to the template $t$ will be denoted by $\z(t)$, see Figure \ref{infzig0}.

\begin{figure}[h]
    \centering
    \scalebox{0.7}{

\tikzset{every picture/.style={line width=0.75pt}} 

\begin{tikzpicture}[x=0.75pt,y=0.75pt,yscale=-1,xscale=1]

\draw  [draw opacity=0] (81.5,94) -- (171.5,94) -- (171.5,104) -- (81.5,104) -- cycle ; \draw   (91.5,94) -- (91.5,104)(101.5,94) -- (101.5,104)(111.5,94) -- (111.5,104)(121.5,94) -- (121.5,104)(131.5,94) -- (131.5,104)(141.5,94) -- (141.5,104)(151.5,94) -- (151.5,104)(161.5,94) -- (161.5,104) ; \draw    ; \draw   (81.5,94) -- (171.5,94) -- (171.5,104) -- (81.5,104) -- cycle ;
\draw  [draw opacity=0] (171.5,104) -- (181.5,104) -- (181.5,194) -- (171.5,194) -- cycle ; \draw    ; \draw   (171.5,114) -- (181.5,114)(171.5,124) -- (181.5,124)(171.5,134) -- (181.5,134)(171.5,144) -- (181.5,144)(171.5,154) -- (181.5,154)(171.5,164) -- (181.5,164)(171.5,174) -- (181.5,174)(171.5,184) -- (181.5,184) ; \draw   (171.5,104) -- (181.5,104) -- (181.5,194) -- (171.5,194) -- cycle ;
\draw  [draw opacity=0] (171.5,94) -- (181.5,94) -- (181.5,104) -- (171.5,104) -- cycle ; \draw    ; \draw    ; \draw   (171.5,94) -- (181.5,94) -- (181.5,104) -- (171.5,104) -- cycle ;
\draw  [draw opacity=0][fill={rgb, 255:red, 155; green, 155; blue, 155 }  ,fill opacity=1 ] (71.5,94) -- (81.5,94) -- (81.5,104) -- (71.5,104) -- cycle ; \draw    ; \draw    ; \draw   (71.5,94) -- (81.5,94) -- (81.5,104) -- (71.5,104) -- cycle ;
\draw  [draw opacity=0][fill={rgb, 255:red, 155; green, 155; blue, 155 }  ,fill opacity=1 ] (71.5,84) -- (81.5,84) -- (81.5,94) -- (71.5,94) -- cycle ; \draw    ; \draw    ; \draw   (71.5,84) -- (81.5,84) -- (81.5,94) -- (71.5,94) -- cycle ;
\draw  [draw opacity=0][fill={rgb, 255:red, 155; green, 155; blue, 155 }  ,fill opacity=1 ] (171.5,194) -- (181.5,194) -- (181.5,204) -- (171.5,204) -- cycle ; \draw    ; \draw    ; \draw   (171.5,194) -- (181.5,194) -- (181.5,204) -- (171.5,204) -- cycle ;
\draw  [draw opacity=0][fill={rgb, 255:red, 155; green, 155; blue, 155 }  ,fill opacity=1 ] (381.5,224) -- (391.5,224) -- (391.5,244) -- (381.5,244) -- cycle ; \draw    ; \draw   (381.5,234) -- (391.5,234) ; \draw   (381.5,224) -- (391.5,224) -- (391.5,244) -- (381.5,244) -- cycle ;
\draw  [draw opacity=0] (191.5,204) -- (281.5,204) -- (281.5,214) -- (191.5,214) -- cycle ; \draw   (201.5,204) -- (201.5,214)(211.5,204) -- (211.5,214)(221.5,204) -- (221.5,214)(231.5,204) -- (231.5,214)(241.5,204) -- (241.5,214)(251.5,204) -- (251.5,214)(261.5,204) -- (261.5,214)(271.5,204) -- (271.5,214) ; \draw    ; \draw   (191.5,204) -- (281.5,204) -- (281.5,214) -- (191.5,214) -- cycle ;
\draw  [draw opacity=0] (291.5,224) -- (381.5,224) -- (381.5,234) -- (291.5,234) -- cycle ; \draw   (301.5,224) -- (301.5,234)(311.5,224) -- (311.5,234)(321.5,224) -- (321.5,234)(331.5,224) -- (331.5,234)(341.5,224) -- (341.5,234)(351.5,224) -- (351.5,234)(361.5,224) -- (361.5,234)(371.5,224) -- (371.5,234) ; \draw    ; \draw   (291.5,224) -- (381.5,224) -- (381.5,234) -- (291.5,234) -- cycle ;
\draw  [draw opacity=0][fill={rgb, 255:red, 155; green, 155; blue, 155 }  ,fill opacity=1 ] (281.5,214) -- (291.5,214) -- (291.5,224) -- (281.5,224) -- cycle ; \draw    ; \draw    ; \draw   (281.5,214) -- (291.5,214) -- (291.5,224) -- (281.5,224) -- cycle ;
\draw  [draw opacity=0][fill={rgb, 255:red, 155; green, 155; blue, 155 }  ,fill opacity=1 ] (281.5,224) -- (291.5,224) -- (291.5,234) -- (281.5,234) -- cycle ; \draw    ; \draw    ; \draw   (281.5,224) -- (291.5,224) -- (291.5,234) -- (281.5,234) -- cycle ;
\draw  [draw opacity=0][fill={rgb, 255:red, 155; green, 155; blue, 155 }  ,fill opacity=1 ] (181.5,194) -- (191.5,194) -- (191.5,204) -- (181.5,204) -- cycle ; \draw    ; \draw    ; \draw   (181.5,194) -- (191.5,194) -- (191.5,204) -- (181.5,204) -- cycle ;
\draw  [draw opacity=0][fill={rgb, 255:red, 155; green, 155; blue, 155 }  ,fill opacity=1 ] (181.5,204) -- (191.5,204) -- (191.5,214) -- (181.5,214) -- cycle ; \draw    ; \draw    ; \draw   (181.5,204) -- (191.5,204) -- (191.5,214) -- (181.5,214) -- cycle ;
\draw  [draw opacity=0][fill={rgb, 255:red, 155; green, 155; blue, 155 }  ,fill opacity=1 ] (281.5,204) -- (291.5,204) -- (291.5,214) -- (281.5,214) -- cycle ; \draw    ; \draw    ; \draw   (281.5,204) -- (291.5,204) -- (291.5,214) -- (281.5,214) -- cycle ;
\draw  [draw opacity=0][fill={rgb, 255:red, 155; green, 155; blue, 155 }  ,fill opacity=1 ] (391.5,234) -- (401.5,234) -- (401.5,254) -- (391.5,254) -- cycle ; \draw    ; \draw   (391.5,244) -- (401.5,244) ; \draw   (391.5,234) -- (401.5,234) -- (401.5,254) -- (391.5,254) -- cycle ;
\draw  [draw opacity=0] (491.5,254) -- (501.5,254) -- (501.5,264) -- (491.5,264) -- cycle ; \draw    ; \draw    ; \draw   (491.5,254) -- (501.5,254) -- (501.5,264) -- (491.5,264) -- cycle ;
\draw  [draw opacity=0] (401.5,254) -- (491.5,254) -- (491.5,264) -- (401.5,264) -- cycle ; \draw   (411.5,254) -- (411.5,264)(421.5,254) -- (421.5,264)(431.5,254) -- (431.5,264)(441.5,254) -- (441.5,264)(451.5,254) -- (451.5,264)(461.5,254) -- (461.5,264)(471.5,254) -- (471.5,264)(481.5,254) -- (481.5,264) ; \draw    ; \draw   (401.5,254) -- (491.5,254) -- (491.5,264) -- (401.5,264) -- cycle ;
\draw  [draw opacity=0] (491.5,264) -- (501.5,264) -- (501.5,354) -- (491.5,354) -- cycle ; \draw    ; \draw   (491.5,274) -- (501.5,274)(491.5,284) -- (501.5,284)(491.5,294) -- (501.5,294)(491.5,304) -- (501.5,304)(491.5,314) -- (501.5,314)(491.5,324) -- (501.5,324)(491.5,334) -- (501.5,334)(491.5,344) -- (501.5,344) ; \draw   (491.5,264) -- (501.5,264) -- (501.5,354) -- (491.5,354) -- cycle ;
\draw  [draw opacity=0][fill={rgb, 255:red, 155; green, 155; blue, 155 }  ,fill opacity=1 ] (391.5,254) -- (401.5,254) -- (401.5,264) -- (391.5,264) -- cycle ; \draw    ; \draw    ; \draw   (391.5,254) -- (401.5,254) -- (401.5,264) -- (391.5,264) -- cycle ;
\draw  [draw opacity=0][fill={rgb, 255:red, 155; green, 155; blue, 155 }  ,fill opacity=1 ] (491.5,354) -- (511.5,354) -- (511.5,364) -- (491.5,364) -- cycle ; \draw   (501.5,354) -- (501.5,364) ; \draw    ; \draw   (491.5,354) -- (511.5,354) -- (511.5,364) -- (491.5,364) -- cycle ;
\draw  [draw opacity=0] (501.5,364) -- (511.5,364) -- (511.5,454) -- (501.5,454) -- cycle ; \draw    ; \draw   (501.5,374) -- (511.5,374)(501.5,384) -- (511.5,384)(501.5,394) -- (511.5,394)(501.5,404) -- (511.5,404)(501.5,414) -- (511.5,414)(501.5,424) -- (511.5,424)(501.5,434) -- (511.5,434)(501.5,444) -- (511.5,444) ; \draw   (501.5,364) -- (511.5,364) -- (511.5,454) -- (501.5,454) -- cycle ;

\end{tikzpicture}
}
    \caption{The infinite zigzag $\z(t)$ for the template $\\$ $\\$ \phantom{FIGURE. 2.\,\,The infinite z}$t=\m{1}\pp\mm\p{1}\m{1}\pp\m{2}\pp\m{1}\p{1}\m{2}\pp\mm\p{1}\mm.\\$
    $\\$
    White strips of boxes represent infinitely long rows and columns while grey boxes represent zigzags corresponding to finite clusters of $t$.}
    \label{infzig0}
\end{figure} 

To every template $t$ we associate a coideal $\zig(t)$ of the zigzag graph, which is by definition of the form $\zig_{\tau}$ for some infinite path $\tau$, see the paragraph above Proposition \ref{co-ideal}. In order to define this path $\tau$, we replace infinite rows and columns in the infinite zigzag $\z(t)$ with long enough (but finite) rows and columns. So, we obtain a sequence of increasing zigzags. Then $\tau$ is any path in the zigzag graph that goes through all these zigzags. Equivalently, $\tau$ is any path that goes through all zigzags corresponding to the binary words obtained from $t$ by replacing infinite clusters with long enough (but finite) blocks. Any such path $\tau$ completely "fills" the infinite zigzag $\z(t)$ that is, starting from some point, $\tau$ looks like a tuple of rows and columns, some of which grow infinitely large while others stay frozen; the frozen rows and columns correspond to finite clusters of $t$. Note that coideals $\zig(t_1)$ and $\zig(t_2)$ coincide if and only if templates $t_1$ and $t_2$ coincide. Below it will be useful sometimes to identify a template $t$ with the corresponding coideal $\zig(t)$. Moreover, it will be convenient to view an infinite zigzag as a sequence of growing finite zigzags.

\begin{proposition}\label{co-ideals of the zigzag graph}
Each proper primitive saturated coideal of the zigzag graph is of the form $\zig(t)$ for some template $t$, which is uniquely defined.
\end{proposition}
\begin{proof}
Let $J$ be a proper saturated primitive coideal of the zigzag graph. From Proposition \ref{co-ideal} it follows that there exists a path $\tau$ such that $J=\zig_{\tau}$. Since $J$ is proper, it follows that the number of blocks in binary words corresponding to zigzags from the path $\tau$ is uniformly bounded along the path. Then we form a template $t$ in the following way. Bounded blocks of binary words from $\tau$ correspond to finite clusters of $t$, and unbounded blocks correspond to infinite clusters.
\end{proof}

\subsubsection{Templates and ideals of $\qsym$}\label{parparpar}
For a template $t$ we consider the following linear subspace $I_t=\spn{\R}{F_{\la}\mid \la\notin\zig(t)}$ of $\qsym$. Note that $I_t$ is a graded ideal of $\qsym$, due to $F_{\la}F_{\mu}\in\spn{\R}{F_{\nu}\mid\nu\geq \la,\mu}$, see. \cites[p. 35, (3.13)]{quaisymmetric_book}.     
\begin{observation}\label{observ fin type2}
Let $t_1,\ldots,t_k$ be some templates. Then $\zig\lr{t_1}\times\ldots\times\zig\lr{t_k}$ is a multiplicative graph, which algebra is $\spn{\R}{F_{\la(1)}\otimes\ldots\otimes F_{\la(k)}\mid\la(i)\in\zig\lr{t_i}}=\quot{\qsym}{I_{t_1}}\otimes\ldots\otimes\quot{\qsym}{I_{t_k}}$, see Definition \ref{direct product graph def}. 
\end{observation}
    
\subsection{Zero sets of finite harmonic functions}
    
From Proposition \ref{prop prim coideal} it follows that the support of a finite indecomposable harmonic function on the zigzag graph is a primitive coideal. If this coideal is proper, then by Proposition \ref{co-ideals of the zigzag graph} it corresponds to a template. Below we specify which finite indecomposable harmonic functions have non-empty zero sets and for them we explicitly describe the corresponding templates.

\subsubsection{Kerov's construction}\label{par fin char}
Recall the definition of finitary oriented paintbox, see Definition 5.2 from \cite{gnedin_olsh2006}.

\begin{definition}\label{space1}
A \textit{finitary oriented paintbox} is a pair $(w_+,w_-)$ of disjoint open subsets of the unit interval $(0,1)$, each comprised of finitely many subintervals and such that the total Lebesgue measure of $w_+$ and $w_-$ equals $1$. The symbol $W_0$ stands for the set of all such pairs.
\end{definition}

Lengths of intervals in $w=(w_+,w_-)\in W_0$ will be denoted by $w_i$. We agree that the intervals are ordered from left to right and $w_1$ denotes the length of the leftmost interval. We say that an interval of $w$ is \textit{positively oriented} if it belongs to $w_+$ and we say that an interval is \textit{negatively oriented} if it belongs to $w_-$.

Kerov's construction produces a finite indecomposable harmonic function $\varphi_w$ on the zigzag graph out of any finitary oriented paintbox $w\in W_0$, see \cite[p. 13-18]{gnedin_olsh2006}. Let us briefly recall this procedure. 

For any zigzag $\la$ we set by definition $$\varphi_w(\la)=F_{\la}(w),$$ 
where
\begin{equation}\label{def fin char}
    F_{\la}(w)=\lr{\psi_1\otimes\ldots\otimes \psi_m}\circ \lr{r_{w_1}\otimes \ldots \otimes r_{w_m}}\circ \Delta^{(m)}\lr{F_{\la}},\ \la\in\zig
\end{equation}
and
\begin{itemize}
    \item $m$ is the number of intervals in $w$;
        
    \item $\Delta^{(m)}\colon \qsym\rightarrow \qsym^{\otimes m}$ is the $m$-th iteration of the comultiplication $\Delta$ in $\qsym$;

    \item $r_t$ is the automorphism of the graded algebra $\qsym$, defined by $r_t\lr{F_{\la}}=t^{|\la|}F_{\la}$;

    \item $\psi_i=\psi_+$, if $w_i$ is positively oriented and $\psi_i=\psi_-$, if $w_i$ is negatively oriented, where  
\begin{equation}\label{psi}
    \psi_+(F_{\la})=\begin{cases}
    1,\ \text{if $\la$ is a row},\\
    
    0\ \text{otherwise},
    \end{cases}\ \ \psi_-(F_{\la})=\begin{cases}
    1\ ,\ \text{if $\la$ is a column},\\
    
    0\ \text{otherwise}.
    \end{cases}
    \end{equation}
\end{itemize}

    Let us consider all splittings of the zigzag $\la$ into $m$ zigzags $\la(1),\ldots,\la(m)$ such that $\la(i)$ is a row, if the interval $w_i$ is positively oriented, and $\la(i)$ is a column, if $w_i$ is negatively oriented. Note that some of these $\la(i)$ may be empty. 
    
    \begin{proposition}\cite[Proposition 5.3]{gnedin_olsh2006}\label{prop prop}
    The following equality holds 
    $$F_{\la}(w)=\summ w_1^{|\la(1)|}w_2^{|\la(2)|}\ldots w_m^{|\la(m)|},$$ 
    where the sum is taken over all splittings of $\la$ mentioned above.
    \end{proposition}
     
Let us denote by $W$ the set of pairs of disjoint open subsets of the unit interval. Note that $W_0$ is a subset of $W$.
     
\begin{theorem}\cite[Theorem 7.5]{gnedin_olsh2006}\label{theorem fin charr}
There is a bijective correspondence $w\mapsto \varphi_w$ between elements of $w\in W$ and indecomposable finite harmonic functions on the zigzag graph. For finitary oriented paintboxes this correspondence is defined by Kerov's construction.
\end{theorem}
    
For any finitary oriented paintbox $w\in W_0$ we denote by $t_w$ the template obtained from $w$ by replacing positively and negatively oriented intervals with symbols $\pp$ and $\mm$ respectively and inserting between any two neighbor infinite symbols of the same type a symbol of the opposite type.
    
    \begin{proposition}\label{kernels}
    Let $w\in W$. Then $\supp(\varphi_w)=\begin{cases}
        \zig(t_w),\ \text{if}\ w\in W_0,\\
        \zig,\ \text{otherwise.}
        \end{cases}$
    \end{proposition}
    \begin{proof}
    Suppose that the finitary oriented paintbox $w$ consists of $m$ intervals. Then Proposition \ref{prop prop} implies that $\varphi_w(\la)>0$ if and only if $\la$ can be represented as a consecutive union of $m$ rows and columns taken in the order proposed by the orientations of intervals of $w$. Thus, $\varphi_w(\la)>0$ if and only if $\la\in\zig(t_w)$. 
    
    Now let $w\in W\backslash W_0$. It suffices to show that $\varphi_w(\la_{2n})>0$ for any $n$, where
    $$\bw(\la_{2n})=\underbrace{+-\ldots+-}_{2n}.$$ 
    For that we will use the oriented paintbox construction from \cite{gnedin_olsh2006}, see Definition 5.4 and the paragraph above Proposition 6.3 in that paper. Following its notation, it remains to prove that the probability $\Pas\lr{\Pi_{2n+1}=\pi_{2n+1}}$ is non-zero, where in one-line notation the permutation $\pi_{2n+1}\in S_{2n+1}$ is given by $$\pi_{2n+1}=1,2n+1,2,2n,3,2n-1,\ldots n-1,n+3,n,n+2,n+1.$$ 
    This fact immediately follows from the next observation. If $w\in W\backslash W_0$ contains infinitely many intervals, then we can place random points inside different intervals in the desired order. But if $w\in W\backslash W_0$ consists of finitely many intervals, then their common length is strictly less then $1$ and we can place our random points inside that complement to $w$, which length is non-zero.
	\end{proof}
	
	\begin{remark}\label{rem zero temp}
	If $w\in W_0$, then the template $t_w$ does not contain finite clusters except those one-symbol clusters which are not outermost and whose two neighbors are infinite clusters of the same sign. Such templates will be called \textit{finite}. A template which is not finite will be called \textit{semifinite}, see Figure \ref{infzig_pair}.
	\end{remark}
	
\begin{figure}[h]
    \centering
    \import{Pictures}{Picture_examp5.tex}
    \caption{a) $\z(t)$ for the finite template $t=\pp\mm\pp\m{1}\pp\m{1}\pp\mm\p{1}\mm.\\$
    b) $\z(t)$ for the semifinite template $t=\m{1}\pp\mm\p{1}\m{1}\pp\m{2}\pp\m{1}\p{1}\m{2}\pp\mm\p{1}\mm.$}
    \label{infzig_pair}
\end{figure} 

\begin{remark}\label{rm2234432}
Proposition \ref{kernels} allows us to think of a finitary oriented paintbox $w=(w_+,w_-)$ as the infinite zigzag $\z(t_w)$ endowed with a tuple of real numbers. Namely, we attach lengths of the intervals from $w_+$ and $w_-$ to infinite rows and columns of $\z(t_w)$ respectively. Moreover, we may identify the infinite zigzag $\z(t_w)$ with an infinite path $\tau$ which completely "fills" this zigzag, see the paragraph above Proposition \ref{co-ideals of the zigzag graph}. Starting from some point, this path $\tau$ looks like a collection of growing rows and columns, hence we can treat the lengths of intervals from $w$ as frequencies of appearing new boxes in that rows and columns which grow infinitely large.
\end{remark}

\subsubsection{A useful lemma}\label{useful lemma}

Now we would like to discuss a lemma, which we will use to prove a semifinite analog of the Vershik-Kerov ring theorem for indecomposable semifinite harmonic functions on the zigzag graph, Proposition \ref{prop mult for zig}.

Let $u$ be an $m$-tuple of adjacent oriented intervals; their lengths will be denoted by $u_1,\ldots,u_m$. The only thing that differs $u$ from a finitary oriented paintbox is that we do not impose any restrictions on lengths of the intervals. For any zigzag $\la$ the expression $F_{\la}(u)$ is defined by the formula from Proposition \ref{prop prop}. Equivalently, we can define this expression by Kerov's construction \eqref{def fin char}. Then it is obvious that $F_{\la}\mapsto F_{\la}(u)$ is a homomorphism of algebras $\qsym\rightarrow \R$. We can also define a template $t_u$ in the same manner as for finitary oriented paintboxes, see the paragraph above Proposition \ref{kernels}. 

Now suppose that $\la\in\zig(t_u)$ and $\bw(\la)$ contains as many blocks as possible. Then each block of $\bw(\la)$ either corresponds to an interval of $u$ or it is a one-symbol block that is placed between two blocks corresponding to equally oriented intervals. The blocks corresponding to intervals of $u$ will be denoted by $\La_1,\ldots,\La_m$. Recall that $|\La_i|$ denotes the number of symbols in the binary word $\La_i$. 

Let us introduce the following notation $n_1\lr{\la}=|\La_1|$, $n_m\lr{\la}=|\La_m|$ and for any $i=2,\ldots,m-1$ we set  
\begin{equation}
    n_i(\la)=\begin{cases}
|\La_i|+1,\ \text{if}\ \text{the}\ i-1\text{-st},i\text{-th},i+1\text{-st}\ \text{intervals of}\ u\ \text{are equally oriented},\\
|\La_i|,\ \text{if}\ \text{the}\ i\text{-th}\ \text{interval of}\ u\ \text{has exactly one neighbor of the same orientation},\\
|\La_i|-1\ \text{otherwise}.\\
\end{cases}
\end{equation}
Let us introduce more notation: $s(u)=|S(u)|$, where
$$S(u)=\{i\mid 1\leq i\leq m\colon \text{the}\ i\text{-th}\ \text{and}\ \text{the}\ i+1\text{-st}\ \text{intervals of}\ u\ \text{have different orientations}\}.$$

\begin{lemma}\label{lemma for multiplic}
Assume that $\la\in\zig(t_u)$ and $\bw(\la)$ contains as many blocks as possible. Then 
$$F_{\la}(u)=u_1^{n_1(\la)}\ldots u_m^{n_m(\la)}\cdot\summ_{\rho\in\{0,1\}^{s(u)}}\prodd_{i\in S(u)}u_{i+\rho(i)},$$ 
where the sum is taken over all $s(u)$-tuples $\rho$, consisting of $0$'s and $1$'s.
\end{lemma}
\begin{proof}
The claim follows from the very definition of $F_{\la}(u)$. Namely, the sum corresponds to all possible splittings of $\la$ mentioned above Proposition \ref{prop prop}.
\end{proof}

\subsection{Semifinite templates}\label{section semifin temp}
From Proposition \ref{cunning Boyer} and Observation \ref{observ fin type2} with $k=1$ it follows that for a finite template $t$ the graph $\zig(t)$ possess no strictly positive indecomposable semifinite harmonic functions, hence $\zig(t)$ can not be realised as the support of an indecomposable semifinite harmonic function on the zigzag graph. Thus, below we will be interested only in semifinite templates.

It turns out that for any semifinite template $t$ the coideal $\zig(t)$ can be realised as the support of an indecomposable semifinite harmonic function on the zigzag graph. Moreover, for indecomposable semifinite harmonic functions with the common support $\zig(t)$ the finiteness ideal depends only on $t$. In the present section we describe this finiteness ideal. Some examples are given in the next section.

\begin{definition}
    Let $t$ be a semifinite template. By a \textit{separating cluster} of $t$ we mean a one-symbol cluster which is not an outermost cluster of $t$ and whose two neighbors are infinite clusters of the same sign. By the \textit{zigzag flange} of $t$ we call a tuple of binary words each of which consists of finite but not separating clusters of $t$ standing nearby. The zigzag flange will be denoted by $\fl(t)$.
\end{definition}
For instance, if we take $t=\m{1}\pp\mm\p{1}\m{1}\pp\m{2}\pp\m{1}\p{1}\m{2}\pp\mm\p{1}\mm$, which is the semifinite template from Figure \ref{infzig_pair}, then $\fl(t)=(\m{1},\p{1}\m{1},\m{2},\m{1}\p{1}\m{2})$. So, the words from this zigzag flange correspond to the first, second, third, and fourth grey zigzags on the Figure \ref{infzig_pair}b):

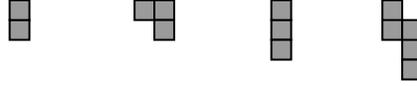
\begin{figure}[h]
    \centering
    \tikzset{every picture/.style={line width=0.75pt}} 

\begin{tikzpicture}[x=0.75pt,y=0.75pt,yscale=-1,xscale=1]

\draw  [draw opacity=0][fill={rgb, 255:red, 155; green, 155; blue, 155 }  ,fill opacity=1 ] (84.5,114) -- (94.5,114) -- (94.5,124) -- (84.5,124) -- cycle ; \draw    ; \draw    ; \draw   (84.5,114) -- (94.5,114) -- (94.5,124) -- (84.5,124) -- cycle ;
\draw  [draw opacity=0][fill={rgb, 255:red, 155; green, 155; blue, 155 }  ,fill opacity=1 ] (84.5,104) -- (94.5,104) -- (94.5,114) -- (84.5,114) -- cycle ; \draw    ; \draw    ; \draw   (84.5,104) -- (94.5,104) -- (94.5,114) -- (84.5,114) -- cycle ;

\draw  [draw opacity=0][fill={rgb, 255:red, 155; green, 155; blue, 155 }  ,fill opacity=1 ] (147.5,104) -- (157.5,104) -- (157.5,114) -- (147.5,114) -- cycle ; \draw    ; \draw    ; \draw   (147.5,104) -- (157.5,104) -- (157.5,114) -- (147.5,114) -- cycle ;
\draw  [draw opacity=0][fill={rgb, 255:red, 155; green, 155; blue, 155 }  ,fill opacity=1 ] (157.5,104) -- (167.5,104) -- (167.5,114) -- (157.5,114) -- cycle ; \draw    ; \draw    ; \draw   (157.5,104) -- (167.5,104) -- (167.5,114) -- (157.5,114) -- cycle ;
\draw  [draw opacity=0][fill={rgb, 255:red, 155; green, 155; blue, 155 }  ,fill opacity=1 ] (157.5,114) -- (167.5,114) -- (167.5,124) -- (157.5,124) -- cycle ; \draw    ; \draw    ; \draw   (157.5,114) -- (167.5,114) -- (167.5,124) -- (157.5,124) -- cycle ;

\draw  [draw opacity=0][fill={rgb, 255:red, 155; green, 155; blue, 155 }  ,fill opacity=1 ] (272.5,104) -- (282.5,104) -- (282.5,124) -- (272.5,124) -- cycle ; \draw    ; \draw   (272.5,114) -- (282.5,114) ; \draw   (272.5,104) -- (282.5,104) -- (282.5,124) -- (272.5,124) -- cycle ;
\draw  [draw opacity=0][fill={rgb, 255:red, 155; green, 155; blue, 155 }  ,fill opacity=1 ] (282.5,114) -- (292.5,114) -- (292.5,134) -- (282.5,134) -- cycle ; \draw    ; \draw   (282.5,124) -- (292.5,124) ; \draw   (282.5,114) -- (292.5,114) -- (292.5,134) -- (282.5,134) -- cycle ;
\draw  [draw opacity=0][fill={rgb, 255:red, 155; green, 155; blue, 155 }  ,fill opacity=1 ] (282.5,134) -- (292.5,134) -- (292.5,144) -- (282.5,144) -- cycle ; \draw    ; \draw    ; \draw   (282.5,134) -- (292.5,134) -- (292.5,144) -- (282.5,144) -- cycle ;

\draw  [draw opacity=0][fill={rgb, 255:red, 155; green, 155; blue, 155 }  ,fill opacity=1 ] (216.5,114) -- (226.5,114) -- (226.5,124) -- (216.5,124) -- cycle ; \draw    ; \draw    ; \draw   (216.5,114) -- (226.5,114) -- (226.5,124) -- (216.5,124) -- cycle ;
\draw  [draw opacity=0][fill={rgb, 255:red, 155; green, 155; blue, 155 }  ,fill opacity=1 ] (216.5,124) -- (226.5,124) -- (226.5,134) -- (216.5,134) -- cycle ; \draw    ; \draw    ; \draw   (216.5,124) -- (226.5,124) -- (226.5,134) -- (216.5,134) -- cycle ;
\draw  [draw opacity=0][fill={rgb, 255:red, 155; green, 155; blue, 155 }  ,fill opacity=1 ] (216.5,104) -- (226.5,104) -- (226.5,114) -- (216.5,114) -- cycle ; \draw    ; \draw    ; \draw   (216.5,104) -- (226.5,104) -- (226.5,114) -- (216.5,114) -- cycle ;

\end{tikzpicture}
\end{figure} 

\begin{definition}
Let $t$ be a semifinite template. By a \textit{section} of $t$ we mean a maximal collection of consecutive clusters that form a finite template.     
\end{definition}
Note that the words from the zigzag flange of $t$ split $t$ into sections.

For the above $t$ the sections are $$t_1=\pp\mm,\ t_2=\pp,\ t_3=\pp,\ t_4=\pp\mm\p{1}\mm$$
and the splitting of $t$ into sections is given by 

$$\lr{\m{1},\ t_1,\ \p{1}\m{1},\ t_2,\ \m{2},\ t_3,\ \m{1}\p{1}\m{2},\ t_4}.$$

In terms of infinite zigzags, we split $\z(t)$ into infinite zigzags corresponding to sections of $t$. For example, $\z(t)$ from Figure \ref{infzig_pair}b) is split by the first, second, third, and fourth grey zigzags into 
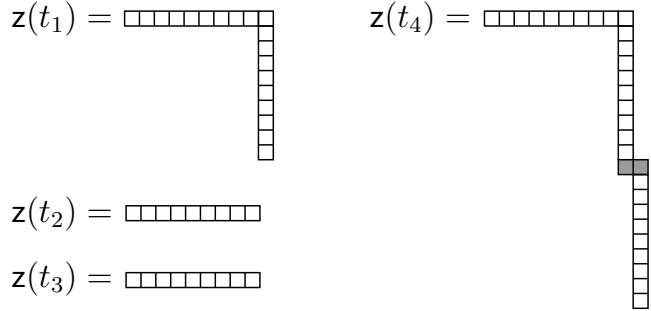
\begin{figure}[h]
    \centering
    \scalebox{0.75}{

\tikzset{every picture/.style={line width=0.75pt}} 

\begin{tikzpicture}[x=0.75pt,y=0.75pt,yscale=-1,xscale=1]

\draw  [draw opacity=0] (59.5,94) -- (149.5,94) -- (149.5,104) -- (59.5,104) -- cycle ; \draw   (69.5,94) -- (69.5,104)(79.5,94) -- (79.5,104)(89.5,94) -- (89.5,104)(99.5,94) -- (99.5,104)(109.5,94) -- (109.5,104)(119.5,94) -- (119.5,104)(129.5,94) -- (129.5,104)(139.5,94) -- (139.5,104) ; \draw    ; \draw   (59.5,94) -- (149.5,94) -- (149.5,104) -- (59.5,104) -- cycle ;
\draw  [draw opacity=0] (149.5,104) -- (159.5,104) -- (159.5,194) -- (149.5,194) -- cycle ; \draw    ; \draw   (149.5,114) -- (159.5,114)(149.5,124) -- (159.5,124)(149.5,134) -- (159.5,134)(149.5,144) -- (159.5,144)(149.5,154) -- (159.5,154)(149.5,164) -- (159.5,164)(149.5,174) -- (159.5,174)(149.5,184) -- (159.5,184) ; \draw   (149.5,104) -- (159.5,104) -- (159.5,194) -- (149.5,194) -- cycle ;
\draw  [draw opacity=0] (149.5,94) -- (159.5,94) -- (159.5,104) -- (149.5,104) -- cycle ; \draw    ; \draw    ; \draw   (149.5,94) -- (159.5,94) -- (159.5,104) -- (149.5,104) -- cycle ;
\draw  [draw opacity=0] (60.5,225) -- (150.5,225) -- (150.5,235) -- (60.5,235) -- cycle ; \draw   (70.5,225) -- (70.5,235)(80.5,225) -- (80.5,235)(90.5,225) -- (90.5,235)(100.5,225) -- (100.5,235)(110.5,225) -- (110.5,235)(120.5,225) -- (120.5,235)(130.5,225) -- (130.5,235)(140.5,225) -- (140.5,235) ; \draw    ; \draw   (60.5,225) -- (150.5,225) -- (150.5,235) -- (60.5,235) -- cycle ;
\draw  [draw opacity=0] (60.5,270) -- (150.5,270) -- (150.5,280) -- (60.5,280) -- cycle ; \draw   (70.5,270) -- (70.5,280)(80.5,270) -- (80.5,280)(90.5,270) -- (90.5,280)(100.5,270) -- (100.5,280)(110.5,270) -- (110.5,280)(120.5,270) -- (120.5,280)(130.5,270) -- (130.5,280)(140.5,270) -- (140.5,280) ; \draw    ; \draw   (60.5,270) -- (150.5,270) -- (150.5,280) -- (60.5,280) -- cycle ;
\draw  [draw opacity=0] (391.5,94) -- (401.5,94) -- (401.5,104) -- (391.5,104) -- cycle ; \draw    ; \draw    ; \draw   (391.5,94) -- (401.5,94) -- (401.5,104) -- (391.5,104) -- cycle ;
\draw  [draw opacity=0] (301.5,94) -- (391.5,94) -- (391.5,104) -- (301.5,104) -- cycle ; \draw   (311.5,94) -- (311.5,104)(321.5,94) -- (321.5,104)(331.5,94) -- (331.5,104)(341.5,94) -- (341.5,104)(351.5,94) -- (351.5,104)(361.5,94) -- (361.5,104)(371.5,94) -- (371.5,104)(381.5,94) -- (381.5,104) ; \draw    ; \draw   (301.5,94) -- (391.5,94) -- (391.5,104) -- (301.5,104) -- cycle ;
\draw  [draw opacity=0] (391.5,104) -- (401.5,104) -- (401.5,194) -- (391.5,194) -- cycle ; \draw    ; \draw   (391.5,114) -- (401.5,114)(391.5,124) -- (401.5,124)(391.5,134) -- (401.5,134)(391.5,144) -- (401.5,144)(391.5,154) -- (401.5,154)(391.5,164) -- (401.5,164)(391.5,174) -- (401.5,174)(391.5,184) -- (401.5,184) ; \draw   (391.5,104) -- (401.5,104) -- (401.5,194) -- (391.5,194) -- cycle ;
\draw  [draw opacity=0][fill={rgb, 255:red, 155; green, 155; blue, 155 }  ,fill opacity=1 ] (391.5,194) -- (411.5,194) -- (411.5,204) -- (391.5,204) -- cycle ; \draw   (401.5,194) -- (401.5,204) ; \draw    ; \draw   (391.5,194) -- (411.5,194) -- (411.5,204) -- (391.5,204) -- cycle ;
\draw  [draw opacity=0] (401.5,204) -- (411.5,204) -- (411.5,294) -- (401.5,294) -- cycle ; \draw    ; \draw   (401.5,214) -- (411.5,214)(401.5,224) -- (411.5,224)(401.5,234) -- (411.5,234)(401.5,244) -- (411.5,244)(401.5,254) -- (411.5,254)(401.5,264) -- (411.5,264)(401.5,274) -- (411.5,274)(401.5,284) -- (411.5,284) ; \draw   (401.5,204) -- (411.5,204) -- (411.5,294) -- (401.5,294) -- cycle ;

\draw (223,86) node [anchor=north west][inner sep=0.75pt]    {\scalebox{1.43}{$\z( t_{4}) =$}};
\draw (-18,86) node [anchor=north west][inner sep=0.75pt]    {\scalebox{1.43}{$\z( t_{1}) =$}};
\draw (-18,261) node [anchor=north west][inner sep=0.75pt]    {\scalebox{1.43}{$\z( t_{3}) =$}};
\draw (-18,217) node [anchor=north west][inner sep=0.75pt]    {\scalebox{1.43}{$\z( t_{2}) =$}};

\end{tikzpicture}
}
\end{figure} 

\begin{remark}
There are analogs of sections and zigzag flanges for the saturated primitive coideals of the Young graph, see the picture on page 148 in \cite{wassermann1981}. Namely, each saturated primitive coideal of the Young graph looks like a thick infinite hook with a flange consisting of a single Young diagram. In our case of the zigzag graph, sections with a zigzag flange play the role of that infinite hook with a Young diagram.
\end{remark}

\begin{definition}
    Let us set $J(t)=\bigcup_{r}\zig(r)$, where the union is taken over all $r$ obtained from $t$ by removing a single symbol from some cluster corresponding to a block of a binary word from the zigzag flange $\fl(t)$.
\end{definition}
Note that $r$ from the definition above may fail to be a template, since $r$ can contain two neighbor clusters of the same sign. Anyway, the construction of $\zig(r)$ remains unchanged. Namely, $\zig(r)$ is the coideal of $\zig$ corresponding to any path passing through the zigzags corresponding to binary words obtained from $r$ by replacing infinite clusters with long enough blocks. This means that we merge two neighbor clusters of the same sign in $r$ into a bigger cluster by adding their lengths.

The ideal $\zig(t)\backslash J(t)$ of $\zig(t)$ is going to be the finiteness ideal of any strictly positive indecomposable semifinite harmonic function on $\zig(t)$. Recall that these functions are in an obvious bijection with indecomposable semifinite harmonic functions on $\zig$ whose support equals $\zig(t)$. Now we would like to describe the ideal $\zig(t)\backslash J(t)$ in more details. 

Let $t$ be a semifinite template with $k$ sections $t_1,\ldots,t_k$. Assume that $\fl(t)=(a_0,\ldots,a_k)$ and the splitting of $t$ into sections looks like

$$t=(a_0,\ t_1,\ a_1,\ \ldots\ , a_{k-1},\ t_k,\ a_k).$$
If $a_0$ or $a_k$ is the empty binary word, then we should merely ignore it in all what follows.

For binary words $a$ and $b$ we write $a>b$ if and only if the number of symbols in $a$ is greater then the number of symbols in $b$ and $b$ can be obtained from $a$ by removing some symbols. The symbol $\sqcup$ denotes the concatenation of binary words. For instance, $-\sqcup +=-+$, $\m{2}\sqcup \m{3}=\m{5}$ and the empty binary word is the identity for $\sqcup$. 

\begin{lemma}\label{lemma injection}
\phantom{}
\begin{enumerate}[leftmargin=3ex, label=\theenumi)]
    \item If $\la\in\zig(t)\backslash J(t)$, then
    $$\bw(\la)=a_0\sqcup \bw(\la^{(1)}) \sqcup a_1\sqcup \ldots \sqcup \bw(\la^{(k)})\sqcup a_k$$
    for some $\la^{(i)}\in\zig(t_i)$, which are uniquely defined.

    \item The following map 
    $$\zig(t)\backslash J(t)\rightarrow \zig(t_1)\times \ldots \times\zig(t_k),\ \ \ \la\mapsto (\la^{(1)},\ldots, \la^{(k)}),$$
    provided by the first part of the lemma, defines an embedding of graded graphs\footnote{By an embedding $f\colon \Ga_1\rightarrow \Ga_2$ of graded graphs we mean an injective map between the sets of vertices such that for any $\la,\mu\in\Ga_1$ we have $\la\nearrow \mu$ if and only if $f(\la)\nearrow f(\mu)$.}, see Definition \ref{direct product graph def}. Moreover, the image of this embedding is an ideal.
\end{enumerate}
\end{lemma}

\begin{proof}
\begin{enumerate}[leftmargin=3ex, label=\theenumi)]
    \item Obviously, we can write
\begin{align}
    \bw(\zig(t))=\ml{\ml{\{}}b_0\sqcup \bw(\la^{(1)}) \sqcup b_1\sqcup \ldots \sqcup\bw(\la^{(k)})\sqcup b_k&\ \ml{\ml{\mid}}\ b_i\leq a_i, \la^{(i)}\in\zig(t_i)\ml{\ml{\}}},\\
    \bw(J(t))=\ml{\ml{\{}}b_0\sqcup \bw(\la^{(1)}) \sqcup b_1\sqcup \ldots \sqcup\bw(\la^{(k)})\sqcup  b_k&\  \ml{\ml{\mid}}\ b_i\leq a_i\ \text{and}\ \exists j\colon b_j<a_j;\ \la^{(i)}\in\zig(t_i)\ml{\ml{\}}}.
\end{align}
So, it suffices to show that if $\la\in\zig(t)\backslash J(t)$ and 
\begin{equation}\label{eq for proof 3.10}
    \begin{multlined}
    \bw(\la)=a_0\sqcup \bw(\la^{(1)}) \sqcup a_1\sqcup \ldots \sqcup \bw(\la^{(k)})\sqcup a_k=\\
    =a_0\sqcup \bw(\wt{\la^{(1)}}) \sqcup a_1\sqcup \ldots \sqcup \bw(\wt{\la^{(k)}})\sqcup a_k,  
    \end{multlined}
\end{equation}
for some $\la^{(i)},\wt{\la^{(i)}}\in\zig(t_i)$, then $\wt{\la^{(i)}}=\la^{(i)}$.

Let us denote by $m$ the natural number such that $\la^{(1)}=\wt{\la^{(1)}},\ldots,\la^{(m-1)}=\wt{\la^{(m-1)}}$, but $\la^{(m)}\neq\wt{\la^{(m)}}$. Equation \eqref{eq for proof 3.10} yields $|\bw(\la^{(m)})|\neq|\bw(\wt{\la^{(m)}})|$, hence we may assume that $|\bw(\wt{\la^{(m)}})|>|\bw(\la^{(m)})|$. Then we can write 
$$\bw(\wt{\la^{(m)}})=\bw(\la^{(m)})\sqcup \delta\in\bw(\zig(t_m))$$
for some non-empty binary word $\delta$. From equation \eqref{eq for proof 3.10} it follows that the first symbol in $\delta$ is the same as in $a_m$. Thus, the condition $\bw(\la^{(m)})\sqcup \delta\in\bw(\zig(t_m))$ implies that $\la\in J(t)$. This contradiction proves the first part of the lemma.

\item Let us denote the map $\la\mapsto(\la^{(1)},\ldots,\la^{(k)})$ by $f$. Let us show that $\la,\mu\in\zig(t)\backslash J(t)$ are joined by an edge if and only if $f\lr{\la}$ and $f\lr{\mu}$ are joined by an edge. Obviously, if $f(\la) \nearrow f(\mu)$, then $\la\nearrow \mu$. Suppose now that $\la\nearrow\mu$. From the first part of the lemma it follows that
$$\bw(\la)=a_0\sqcup \bw(\la^{(1)}) \sqcup a_1\sqcup \ldots \sqcup \bw(\la^{(k)})\sqcup a_k\in\bw\lr{\zig(t)\backslash J(t)}$$
and 
$$\bw(\mu)=a_0\sqcup \bw(\mu^{(1)}) \sqcup a_1\sqcup \ldots \sqcup \bw(\mu^{(k)})\sqcup a_k\in\bw\lr{\zig(t)\backslash J(t)}$$
for some $\la^{(i)},\mu^{(i)}\in\zig(t_i)$. 

The condition $\la\nearrow\mu$ means that we can obtain $\bw(\la)$ by removing a symbol from $\bw(\mu)$. Therefore, this symbol can be deleted only from some $\bw(\mu^{(i)})$. The thing is that we can not remove that symbol from some $a_i$, because then the result belongs to $\bw(J(t))$, but $\la\notin J(t)$. Thus, $f\lr{\la}$ and $f\lr{\mu}$ are joined by an edge.

Suppose that $\la\in\zig(t)\backslash J(t)$ and $(\mu_1,\ldots,\mu_k)\in\zig(t_1)\times\ldots\times\zig(t_k)$. It is straightforward to check that $f(\la)\nearrow (\mu_1,\ldots, \mu_k)$ yields $(\mu_1,\ldots, \mu_k)\in \Im(f)$. Thus, $\Im(f)$ is an ideal.
\end{enumerate}
\end{proof}

The map provided by Lemma \ref{lemma injection} is not surjective by a trivial reason, since the element 
$$(\diameter,\ldots,\diameter)\in\zig(t_1)\times\ldots\times\zig(t_k),$$
corresponding to the empty zigzag at each factor, can not belong to the image. In fact, this is not the only obstacle for this map to be surjective; see examples below.

\subsection{Examples}
For a binary word $a$ we use the following notation $$\zig(t)^{a}=\{\la\in\zig(t)\mid\bw(\la)\geq a \}.$$

\begin{example}\label{eqwer}
Take $t=\p{1}\mm\pp\m{1}\pp$, then its zigzag flange consists of the binary word $a_0=\p{1}$ and the single section is $t_1=\mm\pp\m{1}\pp$. 

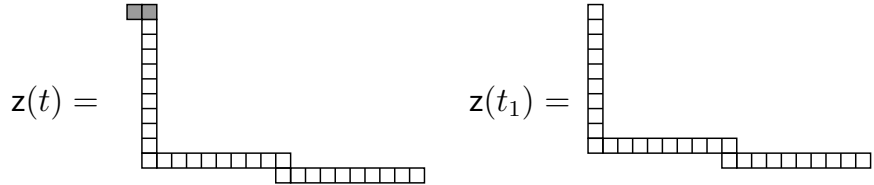
\begin{figure}[H]
    \centering
    \scalebox{0.75}{

\tikzset{every picture/.style={line width=0.75pt}} 

\begin{tikzpicture}[x=0.75pt,y=0.75pt,yscale=-1,xscale=1]

\draw  [draw opacity=0] (91.5,205) -- (181.5,205) -- (181.5,215) -- (91.5,215) -- cycle ; \draw   (101.5,205) -- (101.5,215)(111.5,205) -- (111.5,215)(121.5,205) -- (121.5,215)(131.5,205) -- (131.5,215)(141.5,205) -- (141.5,215)(151.5,205) -- (151.5,215)(161.5,205) -- (161.5,215)(171.5,205) -- (171.5,215) ; \draw    ; \draw   (91.5,205) -- (181.5,205) -- (181.5,215) -- (91.5,215) -- cycle ;
\draw  [draw opacity=0] (81.5,115) -- (91.5,115) -- (91.5,205) -- (81.5,205) -- cycle ; \draw    ; \draw   (81.5,125) -- (91.5,125)(81.5,135) -- (91.5,135)(81.5,145) -- (91.5,145)(81.5,155) -- (91.5,155)(81.5,165) -- (91.5,165)(81.5,175) -- (91.5,175)(81.5,185) -- (91.5,185)(81.5,195) -- (91.5,195) ; \draw   (81.5,115) -- (91.5,115) -- (91.5,205) -- (81.5,205) -- cycle ;
\draw  [draw opacity=0] (81.5,205) -- (91.5,205) -- (91.5,215) -- (81.5,215) -- cycle ; \draw    ; \draw    ; \draw   (81.5,205) -- (91.5,205) -- (91.5,215) -- (81.5,215) -- cycle ;
\draw  [draw opacity=0][fill={rgb, 255:red, 155; green, 155; blue, 155 }  ,fill opacity=1 ] (71.5,105) -- (81.5,105) -- (81.5,115) -- (71.5,115) -- cycle ; \draw    ; \draw    ; \draw   (71.5,105) -- (81.5,105) -- (81.5,115) -- (71.5,115) -- cycle ;
\draw  [draw opacity=0][fill={rgb, 255:red, 155; green, 155; blue, 155 }  ,fill opacity=1 ] (81.5,105) -- (91.5,105) -- (91.5,115) -- (81.5,115) -- cycle ; \draw    ; \draw    ; \draw   (81.5,105) -- (91.5,105) -- (91.5,115) -- (81.5,115) -- cycle ;
\draw  [draw opacity=0] (171.5,215) -- (181.5,215) -- (181.5,225) -- (171.5,225) -- cycle ; \draw    ; \draw    ; \draw   (171.5,215) -- (181.5,215) -- (181.5,225) -- (171.5,225) -- cycle ;
\draw  [draw opacity=0] (181.5,215) -- (271.5,215) -- (271.5,225) -- (181.5,225) -- cycle ; \draw   (191.5,215) -- (191.5,225)(201.5,215) -- (201.5,225)(211.5,215) -- (211.5,225)(221.5,215) -- (221.5,225)(231.5,215) -- (231.5,225)(241.5,215) -- (241.5,225)(251.5,215) -- (251.5,225)(261.5,215) -- (261.5,225) ; \draw    ; \draw   (181.5,215) -- (271.5,215) -- (271.5,225) -- (181.5,225) -- cycle ;
\draw  [draw opacity=0] (391.5,195) -- (481.5,195) -- (481.5,205) -- (391.5,205) -- cycle ; \draw   (401.5,195) -- (401.5,205)(411.5,195) -- (411.5,205)(421.5,195) -- (421.5,205)(431.5,195) -- (431.5,205)(441.5,195) -- (441.5,205)(451.5,195) -- (451.5,205)(461.5,195) -- (461.5,205)(471.5,195) -- (471.5,205) ; \draw    ; \draw   (391.5,195) -- (481.5,195) -- (481.5,205) -- (391.5,205) -- cycle ;
\draw  [draw opacity=0] (381.5,105) -- (391.5,105) -- (391.5,195) -- (381.5,195) -- cycle ; \draw    ; \draw   (381.5,115) -- (391.5,115)(381.5,125) -- (391.5,125)(381.5,135) -- (391.5,135)(381.5,145) -- (391.5,145)(381.5,155) -- (391.5,155)(381.5,165) -- (391.5,165)(381.5,175) -- (391.5,175)(381.5,185) -- (391.5,185) ; \draw   (381.5,105) -- (391.5,105) -- (391.5,195) -- (381.5,195) -- cycle ;
\draw  [draw opacity=0] (381.5,195) -- (391.5,195) -- (391.5,205) -- (381.5,205) -- cycle ; \draw    ; \draw    ; \draw   (381.5,195) -- (391.5,195) -- (391.5,205) -- (381.5,205) -- cycle ;
\draw  [draw opacity=0] (471.5,205) -- (481.5,205) -- (481.5,215) -- (471.5,215) -- cycle ; \draw    ; \draw    ; \draw   (471.5,205) -- (481.5,205) -- (481.5,215) -- (471.5,215) -- cycle ;
\draw  [draw opacity=0] (481.5,205) -- (571.5,205) -- (571.5,215) -- (481.5,215) -- cycle ; \draw   (491.5,205) -- (491.5,215)(501.5,205) -- (501.5,215)(511.5,205) -- (511.5,215)(521.5,205) -- (521.5,215)(531.5,205) -- (531.5,215)(541.5,205) -- (541.5,215)(551.5,205) -- (551.5,215)(561.5,205) -- (561.5,215) ; \draw    ; \draw   (481.5,205) -- (571.5,205) -- (571.5,215) -- (481.5,215) -- cycle ;

\draw (-8,158) node [anchor=north west][inner sep=0.75pt]    {\scalebox{1.43}{$\z( t) =$}};
\draw (300,158) node [anchor=north west][inner sep=0.75pt]    {\scalebox{1.43}{$\z( t_{1}) =$}};

\end{tikzpicture}
}
\end{figure} 

Next, $J(t)=\zig(t_1)$. It is obvious that $J(t)$ and $\zig(t)^{+--}$ do not intersect. Moreover, one can check that $\zig(t)=J(t)\cup \zig(t)^{+--}$. In an expanded form it reads as $$\zig\lr{\p{1}\mm\pp\m{1}\pp}=\zig\lr{\mm\pp\m{1}\pp}\cup \zig\lr{\p{1}\mm\pp\m{1}\pp}^{+--}.$$

Let us prove this. Suppose that $\la\in\zig(t)$. We have to consider several cases.
\begin{itemize}
    \item $\bw(\la)$ has $5$ blocks. Then $\la\in \zig(t)^{+--}$.
    
    \item $\bw(\la)$ has $4$ blocks. Since there are only two types of binary words consisting of $4$ blocks, it follows that either $\bw(\la)=\p{n_1}\m{n_2}\p{n_3}\m{n_4}$ or $\bw(\la)=\m{n_1}\p{n_2}\m{n_3}\p{n_4}$ for some strictly positive integers $n_1,n_2,n_3,n_4$. If the former, then $\la\in \zig(t)^{+--}$; if the latter, then $n_3=1$ and $\la\in\zig(t_1)$.
    
    \item $\bw(\la)$ has $3$ blocks. Then for some strictly positive integers $n_1,n_2,n_3$ one of the following holds:
    \begin{itemize}
        \item $\bw(\la)=\p{n_1}\m{n_2}\p{n_3}$ with $n_2\geq 2$ and then $\la\in \zig(t)^{+--}$;
        
        \item $\bw(\la)=\p{n_1}\m{1}\p{n_3}$ and $\la\in\zig(t_1)$;
        
        \item $\bw(\la)=\m{n_1}\p{n_2}\m{n_3}$ with $n_3\geq 2$ and then $\la\in \zig(t)^{+--}$;
        
        \item $\bw(\la)=\m{n_1}\p{n_2}\m{1}$ and $\la\in\zig(t_1)$.
    \end{itemize}
    
    \item $\bw(\la)$ has $2$ blocks. Then for some strictly positive integers $n_1,n_2$ one of the following holds:
    \begin{itemize}
        \item $\bw(\la)=\m{n_1}\p{n_2}$ and $\la\in\zig(t_1)$;
        
        \item $\bw(\la)=\p{n_1}\m{n_2}$ with $n_2\geq 2$ and then $\la\in \zig(t)^{+--}$;
        
        \item $\bw(\la)=\p{n_1}\m{1}$ and $\la\in\zig(t_1)$;
    \end{itemize}

    \item $\bw(\la)$ consists of a single block. Then $\la\in\zig(t_1)$.
\end{itemize}

Next, we describe the the map provided by Lemma \ref{lemma injection}. 

For any zigzag $\la\in\zig(t)^{+--}$ we can write $\bw(\la)=\p{1}\sqcup \bw(\ol{\la})$ for a unique $\ol{\la}\in \zig(t_1)$ such that $\bw(\ol{\la})$ contains at least two minuses. Thus, the map provided by Lemma \ref{lemma injection} is given by 
$$\zig\lr{\p{1}\mm\pp\m{1}\pp}^{+--}\longrightarrow \zig\lr{\mm\pp\m{1}\pp},\ \ \la\mapsto\ol{\la}.$$

This map is not surjective, since $\bw(\ol{\la})$ contains at least two minuses. The image of this map is the ideal generated by $\m{2}$, that is $\zig\lr{\mm\pp\m{1}\pp}^{--}$.
\end{example}

The ideal $\zig(t)\backslash J(t)$ in the previous example is generated by a single zigzag. The next example shows that this is not the case in general.

\begin{example}\label{ex4263}
Take $t=\m{1}\pp\mm\p{1}\mm\pp\mm\p{1}$, then $a_0=\m{1}$, $a_1=\p{1}$, $t_1=\pp\mm\p{1}\mm\pp\mm$, and
$$J(t)=\zig\lr{\m{1}\pp\mm\p{1}\mm\pp\mm}\cup \zig\lr{\pp\mm\p{1}\mm\pp\mm\p{1}}.$$

\begin{figure}[h]
    \centering
    \import{Pictures}{Picture_examp9.tex}
    \caption{$a)\ \z(t),\ \ b)\ \z\lr{\m{1}\pp\mm\p{1}\mm\pp\mm},\ \ c)\ \z\lr{\pp\mm\p{1}\mm\pp\mm\p{1}}.$}
    \label{fig9}
\end{figure} 

It is not difficult to check that
$$\zig(t)\backslash J(t)=\zig(t)^{-+-+-+-+}\cup \zig(t)^{-\p{2}-\p{2}-+}.$$


\begin{figure}[h]
    \centering
    \tikzset{every picture/.style={line width=0.75pt}} 

\begin{tikzpicture}[x=0.75pt,y=0.75pt,yscale=-1,xscale=1]

\draw  [draw opacity=0] (65.5,60.5) -- (75.5,60.5) -- (75.5,70.5) -- (65.5,70.5) -- cycle ; \draw    ; \draw    ; \draw   (65.5,60.5) -- (75.5,60.5) -- (75.5,70.5) -- (65.5,70.5) -- cycle ;
\draw  [draw opacity=0] (65.5,70.5) -- (75.5,70.5) -- (75.5,80.5) -- (65.5,80.5) -- cycle ; \draw    ; \draw    ; \draw   (65.5,70.5) -- (75.5,70.5) -- (75.5,80.5) -- (65.5,80.5) -- cycle ;
\draw  [draw opacity=0] (75.5,70.5) -- (85.5,70.5) -- (85.5,80.5) -- (75.5,80.5) -- cycle ; \draw    ; \draw    ; \draw   (75.5,70.5) -- (85.5,70.5) -- (85.5,80.5) -- (75.5,80.5) -- cycle ;
\draw  [draw opacity=0] (75.5,80.5) -- (85.5,80.5) -- (85.5,90.5) -- (75.5,90.5) -- cycle ; \draw    ; \draw    ; \draw   (75.5,80.5) -- (85.5,80.5) -- (85.5,90.5) -- (75.5,90.5) -- cycle ;
\draw  [draw opacity=0] (85.5,80.5) -- (95.5,80.5) -- (95.5,90.5) -- (85.5,90.5) -- cycle ; \draw    ; \draw    ; \draw   (85.5,80.5) -- (95.5,80.5) -- (95.5,90.5) -- (85.5,90.5) -- cycle ;
\draw  [draw opacity=0] (85.5,90.5) -- (95.5,90.5) -- (95.5,100.5) -- (85.5,100.5) -- cycle ; \draw    ; \draw    ; \draw   (85.5,90.5) -- (95.5,90.5) -- (95.5,100.5) -- (85.5,100.5) -- cycle ;
\draw  [draw opacity=0] (95.5,90.5) -- (105.5,90.5) -- (105.5,100.5) -- (95.5,100.5) -- cycle ; \draw    ; \draw    ; \draw   (95.5,90.5) -- (105.5,90.5) -- (105.5,100.5) -- (95.5,100.5) -- cycle ;
\draw  [draw opacity=0] (95.5,100.5) -- (105.5,100.5) -- (105.5,110.5) -- (95.5,110.5) -- cycle ; \draw    ; \draw    ; \draw   (95.5,100.5) -- (105.5,100.5) -- (105.5,110.5) -- (95.5,110.5) -- cycle ;
\draw  [draw opacity=0] (105.5,100.5) -- (115.5,100.5) -- (115.5,110.5) -- (105.5,110.5) -- cycle ; \draw    ; \draw    ; \draw   (105.5,100.5) -- (115.5,100.5) -- (115.5,110.5) -- (105.5,110.5) -- cycle ;

\draw  [draw opacity=0] (65.5,147.5) -- (75.5,147.5) -- (75.5,157.5) -- (65.5,157.5) -- cycle ; \draw    ; \draw    ; \draw   (65.5,147.5) -- (75.5,147.5) -- (75.5,157.5) -- (65.5,157.5) -- cycle ;
\draw  [draw opacity=0] (65.5,157.5) -- (75.5,157.5) -- (75.5,167.5) -- (65.5,167.5) -- cycle ; \draw    ; \draw    ; \draw   (65.5,157.5) -- (75.5,157.5) -- (75.5,167.5) -- (65.5,167.5) -- cycle ;
\draw  [draw opacity=0] (75.5,157.5) -- (85.5,157.5) -- (85.5,167.5) -- (75.5,167.5) -- cycle ; \draw    ; \draw    ; \draw   (75.5,157.5) -- (85.5,157.5) -- (85.5,167.5) -- (75.5,167.5) -- cycle ;
\draw  [draw opacity=0] (85.5,157.5) -- (95.5,157.5) -- (95.5,167.5) -- (85.5,167.5) -- cycle ; \draw    ; \draw    ; \draw   (85.5,157.5) -- (95.5,157.5) -- (95.5,167.5) -- (85.5,167.5) -- cycle ;
\draw  [draw opacity=0] (85.5,167.5) -- (95.5,167.5) -- (95.5,177.5) -- (85.5,177.5) -- cycle ; \draw    ; \draw    ; \draw   (85.5,167.5) -- (95.5,167.5) -- (95.5,177.5) -- (85.5,177.5) -- cycle ;
\draw  [draw opacity=0] (95.5,167.5) -- (105.5,167.5) -- (105.5,177.5) -- (95.5,177.5) -- cycle ; \draw    ; \draw    ; \draw   (95.5,167.5) -- (105.5,167.5) -- (105.5,177.5) -- (95.5,177.5) -- cycle ;
\draw  [draw opacity=0] (105.5,167.5) -- (115.5,167.5) -- (115.5,177.5) -- (105.5,177.5) -- cycle ; \draw    ; \draw    ; \draw   (105.5,167.5) -- (115.5,167.5) -- (115.5,177.5) -- (105.5,177.5) -- cycle ;
\draw  [draw opacity=0] (105.5,177.5) -- (115.5,177.5) -- (115.5,187.5) -- (105.5,187.5) -- cycle ; \draw    ; \draw    ; \draw   (105.5,177.5) -- (115.5,177.5) -- (115.5,187.5) -- (105.5,187.5) -- cycle ;
\draw  [draw opacity=0] (115.5,177.5) -- (125.5,177.5) -- (125.5,187.5) -- (115.5,187.5) -- cycle ; \draw    ; \draw    ; \draw   (115.5,177.5) -- (125.5,177.5) -- (125.5,187.5) -- (115.5,187.5) -- cycle ;

\draw (19,73) node [anchor=north west][inner sep=0.75pt]    {$\lambda _{1} =$};
\draw (193,73) node [anchor=north west][inner sep=0.75pt]    {$\bw( \lambda _{1}) =-+-+-+-+$};
\draw (193,154) node [anchor=north west][inner sep=0.75pt]    {$\bw( \lambda _{2}) =-\p{2}-\p{2}-+$};
\draw (19,154) node [anchor=north west][inner sep=0.75pt]    {$\lambda _{2} =$};

\end{tikzpicture}
    \caption{Generators of the ideal $\zig(t)\backslash J(t)$ for $t=\m{1}\pp\mm\p{1}\mm\pp\mm\p{1}$.}
    \label{fig12}
\end{figure} 

Note that the first summand, which is isomorphic to the $5$-th dimensional Pascal graph $\Pas_5$, corresponds to the binary words having the maximal possible number of blocks.

We will show that the ideal $\zig(t)\backslash J(t)$ can not be generated by a single zigzag.

Let us denote by $\downarrow(\la)$ the set of all lower adjacents of $\la$. Then
\begin{equation}
    \begin{multlined}
    \downarrow(-\p{2}-\p{2}-+)=\left\{\p{2}-\p{2}-+,-+-\p{2}-+,-\p{4}-+,\right.\\\left.
    -\p{2}-+-+,-\p{2}-\p{3},-\p{2}-\p{2}-\right\}
    \end{multlined}
\end{equation}
and
\begin{equation}
    \begin{multlined}
    \downarrow(-+-+-+-+)=\left\{+-+-+-+,\m{2}+-+-+,-\p{2}-+-+,-+\m{2}+-+,\right.\\ \left.
    -+-\p{2}-+,-+-+\m{2}+,-+-+-\p{2},-+-+-+-\right\}.
    \end{multlined}
\end{equation}
Thus, 
\begin{equation}
   \downarrow(-+-+-+-+)\ \ml{\cap}  \downarrow(-\p{2}-\p{2}-+)=\left\{-\p{2}-+-+,-+-\p{2}-+\right\}\subset J(t)
\end{equation}
and the ideal $\zig(t)\backslash J(t)$ can not be generated by a single zigzag.

Thus, the graph $\zig(t)\backslash J(t)$ is is not a branching graph but a graded graph, see Figure \ref{fig13}.
\afterpage{
\begin{sidewaysfigure}[h]
    \centering
    \resizebox{1.0\textheight}{!}{
    \tikzset{every picture/.style={line width=0.25pt}}

\begin{tikzpicture}


\node (01) at (4,1){$\ms{\ms{+-+-+-}}$};
\node (02) at (11,1){$\ms{\ms{\overset{2}{+}-\overset{2}{+}-}}$};

\node (11) at (0,5){$\ms{\ms{\overset{2}{+}-+-+-}}$};
\draw (01) to (11);
\draw (02) to (11);

\node (12) at (2,5){$\ms{\ms{+\overset{2}{-}+-+-}}$};
\draw (01) to (12);

\node (13) at (4,5){$\ms{\ms{+-+\overset{2}{-}+-}}$};
\draw (01) to (13);

\node (14) at (6,5){$\ms{\ms{+-+-\overset{2}{+}-}}$};
\draw (01) to (14);
\draw (02) to (14);

\node (15) at (8,5){$\ms{\ms{+-+-+\overset{2}{-}}}$};
\draw (01) to (15);

\node (16) at (11,5){$\ms{\ms{\overset{3}{+}-\overset{2}{+}-}}$};
\draw (02) to (16);
\node (17) at (12.5,5){$\ms{\ms{\overset{2}{+}-\overset{3}{+}-}}$};
\draw (02) to (17);
\node (18) at (14,5){$\ms{\ms{\overset{2}{+}\overset{2}{-}\overset{2}{+}-}}$};
\draw (02) to (18);
\node (19) at (15.5,5){$\ms{\ms{\overset{2}{+}-\overset{2}{+}\overset{2}{-}}}$};
\draw (02) to (19);

\node (21) at (-3,12.5){$\ms{\ms{\overset{3}{+}-+-+-}}$};
\draw (11) to (21);
\draw (16) to (21);

\node (22) at (-1.5,12.5){$\ms{\ms{+\overset{3}{-}+-+-}}$};
\draw (12) to (22);

\node (23) at (0,12.5){$\ms{\ms{+-+\overset{3}{-}+-}}$};
\draw (13) to (23);

\node (24) at (1.5,12.5){$\ms{\ms{+-+-\overset{3}{+}-}}$};
\draw (14) to (24);
\draw (17) to (24);

\node (25) at (3,12.5){$\ms{\ms{+-+-+\overset{3}{-}}}$};
\draw (15) to (25);

\node (26) at (4.5,12.5){$\ms{\ms{\overset{2}{+}\overset{2}{-}+-+-}}$};
\draw (11) to (26);
\draw (12) to (26);
\draw (18) to (26);

\node (27) at (6,12.5){$\ms{\ms{\overset{2}{+}-+\overset{2}{-}+-}}$};
\draw (11) to (27);
\draw (13) to (27);

\node (28) at (7.5,12.5){$\ms{\ms{\overset{2}{+}-+-\overset{2}{+}-}}$};
\draw (11) to (28);
\draw (14) to (28);
\draw (16) to (28);
\draw (17) to (28);
\draw (18) to (28);

\node (29) at (9,12.5){$\ms{\ms{\overset{2}{+}-+-+\overset{2}{-}}}$};
\draw (11) to (29);
\draw (15) to (29);
\draw (19) to (29);

\node (210) at (-1.8,11.5){$\ms{\ms{+\overset{2}{-}+\overset{2}{-}+-}}$};
\draw (12) to (210);
\draw (13) to (210);

\node (211) at (-0.17,11.5){$\ms{\ms{+\overset{2}{-}+-\overset{2}{+}-}}$};
\draw (12) to (211);
\draw (14) to (211);

\node (212) at (1.375,11.5){$\ms{\ms{+\overset{2}{-}+-+\overset{2}{-}}}$};
\draw (12) to (212);
\draw (15) to (212);

\node (213) at (3,11.5){$\ms{\ms{+-+\overset{2}{-}\overset{2}{+}-}}$};
\draw (13) to (213);
\draw (14) to (213);
\draw (18) to (213);

\node (214) at (4.9,11.5){$\ms{\ms{+-+\overset{2}{-}+\overset{2}{-}}}$};
\draw (13) to (214);
\draw (15) to (214);

\node (215) at (6.6,11.5){$\ms{\ms{+-+-\overset{2}{+}\overset{2}{-}}}$};
\draw (14) to (215);
\draw (15) to (215);
\draw (19) to (215);

\node (216) at (13.4,12.5){$\ms{\ms{\overset{4}{+}-\overset{2}{+}-}}$};
\draw (16) to (216);

\node (221) at (14.6,12.5){$\ms{\ms{\overset{2}{+}-\overset{4}{+}-}}$};
\draw (17) to (221);

\node (223) at (15.7,12.5){$\ms{\ms{\overset{2}{+}\overset{3}{-}\overset{2}{+}-}}$};
\draw (18) to (223);

\node (225) at (16.8,12.5){$\ms{\ms{\overset{2}{+}-\overset{2}{+}\overset{3}{-}}}$};
\draw (19) to (225);

\node (217) at (11.4,11.5){$\ms{\ms{\overset{3}{+}\overset{2}{-}\overset{2}{+}-}}$};
\draw (16) to (217);
\draw (18) to (217);

\node (222) at (12.5,11.5){$\ms{\ms{\overset{2}{+}-\overset{3}{+}\overset{2}{-}}}$};
\draw (17) to (222);
\draw (19) to (222);

\node (218) at (13.7,11.5){$\ms{\ms{\overset{3}{+}-\overset{3}{+}-}}$};
\draw (16) to (218);
\draw (17) to (218);

\node (224) at (14.9,11.5){$\ms{\ms{\overset{2}{+}\overset{2}{-}\overset{2}{+}\overset{2}{-}}}$};
\draw (18) to (224);
\draw (19) to (224);

\node (219) at (16.1,11.5){$\ms{\ms{\overset{3}{+}-\overset{2}{+}\overset{2}{-}}}$};
\draw (16) to (219);
\draw (19) to (219);

\node (220) at (17.2,11.5){$\ms{\ms{\overset{2}{+}\overset{2}{-}\overset{3}{+}-}}$};
\draw (17) to (220);
\draw (18) to (220);

\end{tikzpicture}
    }
    \caption{The first three levels of the image of the map $\zig(t)\backslash J(t)\lhook\joinrel\xrightarrow{\ \ \ \ } \zig\lr{\pp\mm\p{1}\mm\pp\mm}$ provided by Lemma \ref{lemma injection} for $t=\m{1}\pp\mm\p{1}\mm\pp\mm\p{1}$, see Example \ref{ex4263}. Zigzags are represented by the corresponding binary words.}
    \label{fig13}
\end{sidewaysfigure}
\clearpage
}

Let us describe the map from Lemma \ref{lemma injection}. 

It is easy to check that the first block of $\bw(\la)$ for $\la\in \zig(t)\backslash J(t)$ must be the negative one-symbol block and the last block of $\bw(\la)$ must be the positive one symbol block. So, the desired map is given by 
$$\zig(t)^{-+-+-+-+}\cup \zig(t)^{-\p{2}-\p{2}-+}\longrightarrow \zig\lr{\pp\mm\p{1}\mm\pp\mm},\ \ \ \la\mapsto\ol{\la},$$
where $t=\m{1}\pp\mm\p{1}\mm\pp\mm\p{1}$ and $\bw(\la)=-\sqcup\bw(\ol{\la})\sqcup +$.
\end{example}

In the next example we propose a sufficient condition for $\zig(t)\backslash J(t)$ to be generated by a single zigzag. This condition is not necessary, because the template $t$ from Example \ref{eqwer} does not satisfy it, but, nevertheless, the corresponding ideal is generated by a single zigzag.

\begin{example}
Let $t$ be a semifinite template. We say that a cluster of $t$ is \textit{internal} if $t$ neither begins nor ends with this cluster. Let us denote by $a_t$ the binary word obtained from $t$ by applying the following rules:
\begin{itemize}
    \item if $t$ begins or ends with an infinite cluster, then this infinite cluster is removed;
    
    \item each internal infinite cluster of $t$ having an infinite neighbour is replaced by a one-symbol cluster of the same sign;
    
    \item each infinite cluster standing between two finite clusters is removed.
\end{itemize}
For instance, $a_t=-\p{2}-+\m{4}$ for $t=\pp\mm\p{2}\mm\pp\m{1}\pp\m{3}.$

It is not difficult to check that $\zig(t)^{a_t}$ is isomorphic as a graded graph to the Pascal graph of an appropriate dimension. 

Suppose that $t$ satisfies the constraints:
\begin{itemize}
    \item $t$ avoids the following patterns $$\pp\mm\p{1}\mm\pp\ \ \ \text{and}\ \ \ \mm\pp\m{1}\pp\mm;$$
    
    \item $t$ does not begin with 
    $$\pp\m{1}\pp\mm\ \ \ \text{or}\ \ \ \mm\p{1}\mm\pp;$$
    
    \item $t$ does not end with 
    $$\mm\pp\m{1}\pp\ \ \ \text{or}\ \ \ \pp\mm\p{1}\mm.$$
\end{itemize}

Then 
$$\zig(t)\backslash J(t)=\zig(t)^{a_t}.$$

Let us prove this. We will analyze the set $\zig(t)\backslash\zig(t)^{a_t}$ in order to prove that it equals $J(t)$. The proof splits into the following parts:
\begin{enumerate}[label=\theenumi)]
    \item if $r$ is obtained from $t$ by removing a separating symbol and merging the two infinite clusters standing near that separating symbol, then there exists $r'$ that can be obtained from $t$ by removing a symbol from a cluster corresponding to a word from the zigzag flange of $t$ and such that $\zig(r)\subset \zig(r')$;
    
    \item $J(t)\subset \zig(t)\backslash\zig(t)^{a_t}$;
    
    \item $J(t)\supset \zig(t)\backslash\zig(t)^{a_t}$
\end{enumerate}

\textit{The first part of the proof.} Since $\zig(r')$ is a coideal of $\zig(t)$, it is sufficient to show that any binary word $a\in \bw(\zig(r))$ which has as many blocks as possible belongs to $\bw(\zig(r'))$. Moreover, we can assume that the blocks of $a$ are large enough that is, the blocks corresponding to finite clusters are of maximal lengths and the blocks corresponding to infinite clusters are of length, say, $N$, where $N$ is large enough. 

We have to deal with the one-symbol cluster of $t$ which is a part of one of the following patterns $\pp\m{1}\pp$ or $\mm\p{1}\mm$. Let us restrict ourselves only to the first case. Then we can rewrite $t$ as
$$t=(p_1, \pp\m{1}\pp, p_2),$$
where $p_1$ and $p_2$ are such that either $p_1$ ends with a finite cluster or $p_2$ starts with it. Note that here we used the constraints on $t$ listed above. The template $r$ looks like
$$r=(p_1, \pp, p_2)$$
and the induced splitting of $a$ reads as
$$a=a^{(1)}\sqcup \p{N}\sqcup a^{(2)},$$
where $a^{(1)}\in\bw(\zig(p_1))$ and $a^{(2)}\in\bw(\zig(p_2))$ contain as many blocks as possible. Then we can obtain the desired template $r'$ by removing a symbol from that finite cluster of $t$ with which $p_1$ ends or $p_2$ begins.  

\textit{The second part of the proof.} Suppose that $a\in \bw(\zig(r))$ is a binary word consisting of maximal possible number of blocks which are large enough, where $r$ is obtained from $t$ by removing a symbol from a cluster corresponding to a word from the zigzag flange of $t$. We will assume that the blocks of $a$ corresponding to infinite clusters of $t$ are of the same length, which we denote by $N$. Then the number of infinite clusters in $r$ is the same as in $t$ and equals the number of blocks of length $N$ in $a$. These blocks of length $N$ split $a$ into parts almost all of which are binary words from the zigzag flange of $t$, except one part, which differs from a binary word from $\fl(t)$ by a single symbol. We may write this splitting and the induced splitting of $a_t$ as
$$a=a^{(1)}\sqcup \beta\sqcup a^{(2)}$$
$$a_t=a_t^{(1)}\sqcup \alpha \sqcup a_t^{(2)},$$
where $\alpha$ is a binary word from $\fl(t)$ and $\beta\nearrow \alpha$. Then one can easily see that $a^{(i)}\geq a_t^{(i)}$ for $i=1,2$, but $$a^{(1)}\not\geq a_t^{(1)}\sqcup \delta_1\ \ \ \ \ \text{and}\ \ \ \ \ a^{(2)}\not\geq \delta_2\sqcup a_t^{(2)},$$ 
where $\delta_1$ and $\delta_2$ denote single symbols from the first and the last clusters of $\alpha$. 

Thus, $a\not\geq a_t$ and the claim follows.

\textit{The third part of the proof.} For notational simplicity, let us assume that the first cluster of $t$ is of sign plus and the total number of clusters in $t$ is even. We denote this number by $n$. Then we can write
$$t=\p{k_1}\m{k_2}\ldots \p{k_{n-1}}\m{k_n},$$
where $k_1,\ldots,k_n$ is the tuple of formal multiplicities of $t$, some of which may be infinite.

Furthermore,
$$\bw(\zig(t))=\{\p{l_1}\m{l_2}\ldots \p{l_{n-1}}\m{l_n}\mid l_i\leq k_i\}.$$

Let us denote by $I\subset \{1,2,\ldots,n\}$ the positions of finite clusters of $t$. Then

\begin{equation}
\bw\lr{\zig(t)^{a_t}}=\left\{\p{l_1}\m{l_2}\ldots \p{l_{n-1}}\m{l_n}\ \middle| \begin{aligned}
&\bullet\ l_i=k_i\ \text{if}\ i\in I;\\
&\bullet\ l_i\geq 1\ \text{if}\ i\notin I,\ i\neq 1,n,\ \text{and}\ i+1\notin I\ \text{or}\ i-1\notin I
\end{aligned}\right\}.    
\end{equation}

So, if $\p{l_1}\m{l_2}\ldots \p{l_{n-1}}\m{l_n}\in \bw(\zig(t)\backslash \zig(t)^{a_t})$, then at least one of the following conditions must hold
\begin{enumerate}[label=\theenumi)]
    \item $l_i\leq k_i-1$ for some $i\in I$
    
    \item $l_i=0$ for some $i$ such that $2\leq i\leq n-1$, $i\notin I$, and $i+1\notin I$ or $i-1\notin I$.
\end{enumerate}

Thus, from the first part of the proof it follows that to prove the desired inclusion it suffices to show that any binary word $a=\p{l_1}\m{l_2}\ldots \p{l_{n-1}}\m{l_n}$ satisfying the second condition above belongs to $\bw(\zig(r))$ for some $r$ which can be obtained from $t$ by removing a symbol from a finite cluster. In order to do so, we rewrite $t$ in the following way 
$$t=s_1\sqcup s \sqcup s_2$$
with $s$ being the maximal template consisting only of infinite clusters and containing the $i$-th infinite cluster of $t$, where $i$ is the number from the second condition above for our binary word $a$. Note that then $s_2$ either begins with a finite cluster or is empty; if the latter, then $s_1$ ends with a finite cluster. Anyway, it is obvious that we can obtain the desired $r$ by removing a single symbol from a cluster of $t$ which is neighbour to $s$ and is an outermost cluster of $s_1$ or $s_2$.
\end{example}

\section{Harmonic functions on $\zig(t)$}\label{main section}
From Proposition \ref{cunning Boyer} and Observation \ref{observ fin type2} with $k=1$ it follows that for a finite template $t$ the graph $\zig(t)$ possess no strictly positive indecomposable semifinite harmonic functions, hence $\zig(t)$ can not be realised as the support of an indecomposable semifinite harmonic function on the zigzag graph. Then Propositions \ref{prop prim coideal} and \ref{co-ideals of the zigzag graph} imply that in order to describe all indecomposable semifinite harmonic functions on the zigzag graph it is sufficient to describe all strictly positive indecomposable semifinite harmonic functions on $\zig(t)$ for any semifinite template $t$. That is what we do in the present section.

Now we would like to introduce some strictly positive functions on the graph $\zig(t)$, see Definition \ref{defdefdeff}. Below we prove that they are pairwise distinct and form an exhaustive list of indecomposable semifinite harmonic functions on $\zig(t)$. 

\begin{definition}\label{framed oriented pb}
    By a \textit{semifinite zigzag growth model} we call a pair $(t,w)$, where $t$ is a semifinite template having $m$ infinite clusters and $w=(w_1,\ldots,w_m)$ is an $m$-tuple of positive real numbers such that $w_1+\ldots+w_m=1$.
\end{definition}
\begin{remark}\label{rem456}
We can assume that these real numbers $w_1,\ldots,w_m$ are assigned to infinite clusters of $t$. Then we can identify a semifinite zigzag growth model $(t,w)$ with the infinite zigzag $\z(t)$ endowed with a tuple of frequencies, see Remark \ref{rm2234432}. Furthermore, we can treat this $w$ as a finitary oriented paintbox the $i$-th interval component of which is of length $w_i$; orientation of this interval component is defined by the sign of the corresponding infinite cluster of $t$: the orientation is positive if the cluster is positive and the orientation is negative if the cluster is negative.  
\end{remark}

Let $t$ have $k$ sections $t_1,\ldots,t_k$. Assume that $\fl(t)=(a_0,\ldots,a_k)$ and the splitting of $t$ into sections looks like

$$t=(a_0,\ t_1,\ a_1,\ \ldots\ , a_{k-1},\ t_k,\ a_k).$$
If $a_0$ or $a_k$ is the empty binary word, then we should merely ignore it in all what follows.

Let $(t,w)$ be a semifinite zigzag growth model. The splitting of $t$ into sections gives us a splitting of $w$
$$w=v_1\sqcup\ldots \sqcup v_k,$$
where each $v_i$ is a tuple of real numbers from $w=(w_1,\ldots,w_m)$ corresponding to the infinite clusters of $t_i$. Note that we may treat each $v_i$ as a collection of oriented subintervals of $(0,1)$; the only thing that differs $v_i$ from a finitary oriented paintbox is the total length of intervals from $v_i$, which may not be equal to $1$.

\begin{definition}\label{defdefdeff}
For any $\la\in\zig(t)$ we set
    $$\varphi_{t,w}(\la)=\begin{cases}
    F_{\la^{(1)}}(v_1)\cdot\ldots\cdot F_{\la^{(k)}}(v_k),\ \text{if}\ \la\in\zig(t)\backslash J(t)\\
    +\infty,\ \text{if}\ \la\in J(t),
    \end{cases}\\$$
    where $\la\mapsto (\la^{(1)},\ldots,\la^{(k)})$ is the map provided by Lemma \ref{lemma injection} and $F_{\la^{(i)}}(v_i)$ is defined by Kerov's construction \eqref{def fin char} or by the formula from Proposition \ref{prop prop}. 
\end{definition}

\begin{example}\label{exexex}
Take $t=\pp\m{1}\p{1}\mm$. Then $a_0$ and $a_2$ are empty binary words, $a_1=\m{1}\p{1}$, $t_1=\pp$, and $t_2=\mm$. 

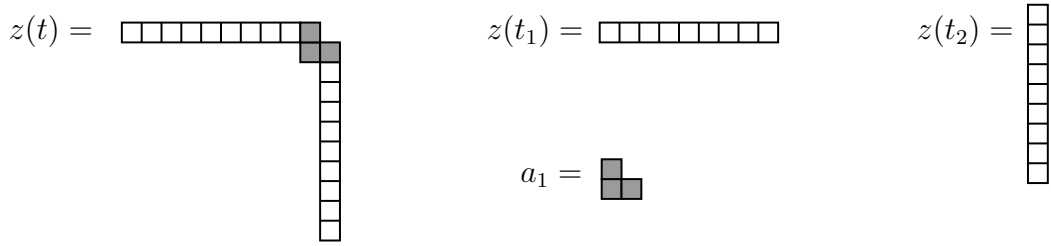
\begin{figure}[h]
    \centering
    \tikzset{every picture/.style={line width=0.75pt}} 

\begin{tikzpicture}[x=0.75pt,y=0.75pt,yscale=-1,xscale=1]

\draw  [draw opacity=0] (99.5,40) -- (189.5,40) -- (189.5,50) -- (99.5,50) -- cycle ; \draw   (109.5,40) -- (109.5,50)(119.5,40) -- (119.5,50)(129.5,40) -- (129.5,50)(139.5,40) -- (139.5,50)(149.5,40) -- (149.5,50)(159.5,40) -- (159.5,50)(169.5,40) -- (169.5,50)(179.5,40) -- (179.5,50) ; \draw    ; \draw   (99.5,40) -- (189.5,40) -- (189.5,50) -- (99.5,50) -- cycle ;
\draw  [draw opacity=0] (199.5,60) -- (209.5,60) -- (209.5,150) -- (199.5,150) -- cycle ; \draw    ; \draw   (199.5,70) -- (209.5,70)(199.5,80) -- (209.5,80)(199.5,90) -- (209.5,90)(199.5,100) -- (209.5,100)(199.5,110) -- (209.5,110)(199.5,120) -- (209.5,120)(199.5,130) -- (209.5,130)(199.5,140) -- (209.5,140) ; \draw   (199.5,60) -- (209.5,60) -- (209.5,150) -- (199.5,150) -- cycle ;
\draw  [draw opacity=0][fill={rgb, 255:red, 155; green, 155; blue, 155 }  ,fill opacity=1 ] (189.5,50) -- (199.5,50) -- (199.5,60) -- (189.5,60) -- cycle ; \draw    ; \draw    ; \draw   (189.5,50) -- (199.5,50) -- (199.5,60) -- (189.5,60) -- cycle ;
\draw  [draw opacity=0][fill={rgb, 255:red, 155; green, 155; blue, 155 }  ,fill opacity=1 ] (189.5,40) -- (199.5,40) -- (199.5,50) -- (189.5,50) -- cycle ; \draw    ; \draw    ; \draw   (189.5,40) -- (199.5,40) -- (199.5,50) -- (189.5,50) -- cycle ;
\draw  [draw opacity=0][fill={rgb, 255:red, 155; green, 155; blue, 155 }  ,fill opacity=1 ] (199.5,50) -- (209.5,50) -- (209.5,60) -- (199.5,60) -- cycle ; \draw    ; \draw    ; \draw   (199.5,50) -- (209.5,50) -- (209.5,60) -- (199.5,60) -- cycle ;

\draw  [draw opacity=0] (340.5,40) -- (430.5,40) -- (430.5,50) -- (340.5,50) -- cycle ; \draw   (350.5,40) -- (350.5,50)(360.5,40) -- (360.5,50)(370.5,40) -- (370.5,50)(380.5,40) -- (380.5,50)(390.5,40) -- (390.5,50)(400.5,40) -- (400.5,50)(410.5,40) -- (410.5,50)(420.5,40) -- (420.5,50) ; \draw    ; \draw   (340.5,40) -- (430.5,40) -- (430.5,50) -- (340.5,50) -- cycle ;
\draw  [draw opacity=0] (556.5,31) -- (566.5,31) -- (566.5,121) -- (556.5,121) -- cycle ; \draw    ; \draw   (556.5,41) -- (566.5,41)(556.5,51) -- (566.5,51)(556.5,61) -- (566.5,61)(556.5,71) -- (566.5,71)(556.5,81) -- (566.5,81)(556.5,91) -- (566.5,91)(556.5,101) -- (566.5,101)(556.5,111) -- (566.5,111) ; \draw   (556.5,31) -- (566.5,31) -- (566.5,121) -- (556.5,121) -- cycle ;

\draw  [draw opacity=0][fill={rgb, 255:red, 155; green, 155; blue, 155 }  ,fill opacity=1 ] (341.5,119) -- (351.5,119) -- (351.5,129) -- (341.5,129) -- cycle ; \draw    ; \draw    ; \draw   (341.5,119) -- (351.5,119) -- (351.5,129) -- (341.5,129) -- cycle ;
\draw  [draw opacity=0][fill={rgb, 255:red, 155; green, 155; blue, 155 }  ,fill opacity=1 ] (341.5,109) -- (351.5,109) -- (351.5,119) -- (341.5,119) -- cycle ; \draw    ; \draw    ; \draw   (341.5,109) -- (351.5,109) -- (351.5,119) -- (341.5,119) -- cycle ;
\draw  [draw opacity=0][fill={rgb, 255:red, 155; green, 155; blue, 155 }  ,fill opacity=1 ] (351.5,119) -- (361.5,119) -- (361.5,129) -- (351.5,129) -- cycle ; \draw    ; \draw    ; \draw   (351.5,119) -- (361.5,119) -- (361.5,129) -- (351.5,129) -- cycle ;

\draw (282,35) node [anchor=north west][inner sep=0.75pt]    {$z( t_{1}) =$};
\draw (41,35) node [anchor=north west][inner sep=0.75pt]    {$z( t) =$};
\draw (499,35) node [anchor=north west][inner sep=0.75pt]    {$z( t_{2}) =$};
\draw (299,112) node [anchor=north west][inner sep=0.75pt]    {$a_{1} =$};

\end{tikzpicture}
\end{figure} 

Next, $J(t)=\zig\lr{\pp\mm}$ and 
$$\zig(t)\backslash J(t)=\zig\lr{\pp\m{1}\p{1}\mm}^{-+}.$$
Recall that the superscript denotes the zigzags which binary words contain $-+$.

The map 
$$\zig\lr{\pp\m{1}\p{1}\mm}^{-+}\longrightarrow \zig\lr{\pp}\times\zig\lr{\mm}$$
provided by Lemma \ref{lemma injection} is given by $\p{n}-+\m{m}\mapsto (\p{n},\m{m})$ and turns out to be as surjective as possible, see the paragraph below the proof of Lemma \ref{lemma injection}.

Let $w_1$ and $w_2$ be real positive numbers such that $w_1+w_2=1$. Then

$$\varphi_{t,w}(\la)=\begin{cases}
w_1^{n+1} w_2^{m+1},\ \text{if}\ \bw(\la)=\p{n}-+\m{m}\ \text{with}\ n,m\geq0,\\
\\
+\infty,\ \text{if}\ \la\in\zig\lr{\pp\mm}. 
\end{cases}$$

It is straightforward to check that $\varphi_{t,w}$ is a harmonic function. Let us show that it is semifinite. This means that for any $\p{n}\m{m}\in \bw(J(t))$ we have to find an approximating sequence, see Remark \ref{remark 321}. Note that we can assume that these $n$ and $m$ are large enough, since if $\la\geq \mu$ and $\{a_N\}_{N\geq 1}$ is an approximating sequence for $\la$, then $\{a_N\}_{N\geq 1}$ is an approximating sequence for $\mu$ as well. In fact, we will use only the bound $n,m\geq 2$. Below we treat binary word belonging to $\bw(\zig(t))$ as elements of $\K(\zig(t))$, see Section \ref{qwer}.

We argue that 
$$a_N=N\cdot \p{n-1}-+\ \m{m-1}$$
for $N\geq 1$ form an approximating sequence for $\p{n}\m{m}$. Let us prove this. 

Since 
$$\varphi_{t,w}(a_N)=N\cdot w_1^n w_2^m$$ 
it follows that $$\varphi_{t,w}(a_N)<+\infty$$
and 
$$\varphi_{t,w}(a_N)\to +\infty\ \text{as}\ N\to+\infty.$$ 

Thus, it suffices to show that $\p{n}\m{m}\geq_K a_N$; for the definition of $\geq_K$ see Section \ref{qwer}.

By the harmonicity condition for any $N\geq 1$ we can write

\begin{equation}
    \begin{multlined}
        \p{n}\m{m}=\summ_{\substack{\la\in\zig(t)\colon\\ |\la|=n+m+N+1}} \dim(\p{n}\m{m},\la)\, \cdot\la\geq_K\\
        \ \ \ \ \ \ \ \ \ \ \ \ \ \ \ \ \ \ \ \ \ \ \ \ \ \ \ \ \ \ \ \ \ \ \summ_{\substack{n_1,m_1\geq 0\colon\\ n_1+m_1=N}}\dim\lr{\p{n}\m{m},\p{n_1+n-1}-+\m{m_1+m-1}}\cdot\p{n_1+n-1}-+\ \m{m_1+m-1}.
    \end{multlined}
\end{equation}

The harmonicity condition also implies
$$a_N=N\cdot \summ_{\substack{n_1,m_1\geq 0\colon\\ n_1+m_1=N}}\dim\lr{\p{n-1}-+\ \m{m-1},\p{n_1+n-1}-+\ \m{m_1+m-1}}\cdot\p{n_1+n-1}-+\ \m{m_1+m-1}.$$

Note that 
$$\dim\lr{\p{n-1}-+\ \m{m-1},\p{n_1+n-1}-+\ \m{m_1+m-1}}=\binom{N}{n_1,m_1},$$
where $\binom{a+b}{a,b}$ denotes the binomial coefficient $\cfrac{(a+b)!}{a!b!}$.

Then the desired claim immediately follows from the next observation.

\begin{equation}\label{tuwer}
    \dim\lr{\p{n}\m{m},\p{n_1+n-1}-+\m{m_1+m-1}}=n_1\cdot\binom{N}{n_1,m_1}+m_1\cdot\binom{N}{n_1,m_1}=N\cdot\binom{N}{n_1,m_1}.
\end{equation}

Let us explain the origin of the first summand in \eqref{tuwer}. It equals the number of the following paths in the zigzag graph going from $\p{n}\m{m}$ to $\p{n_1+n-1}-+\m{m_1+m-1}$. Up to some point we increase only the two blocks of $\p{n}\m{m}$ and then we add the plus that increases the number of blocks in the word; after that, we increase only the outermost blocks again. So, we may treat each such path as an ordered collection of pluses and minuses, which we add to the original word $\p{n}\m{m}$; one plus, which increases the number of blocks, is marked. This collection is of length $N$ and there are $n_1$ pluses in it. Only two things that may vary are positions of ordinary pluses, which we add to the leftmost block, and the position of that special plus that increases the number of blocks. So, we have to choose $n_1$ elements from the set of $N$ elements and mark one of the chosen elements. The number of all such choices equals 
$$n_1\cdot\binom{N}{n_1,m_1}.$$

The second summand in \eqref{tuwer} comes from the similar picture but for minuses. In fact, we have just proved the inequality $\geq$ in \eqref{tuwer}. But one can easily see that this is indeed an equality for $n,m\geq2$.

Finally, we shall prove that $\varphi_{t,w}$ is indecomposable. Suppose that $\psi$ is a finite or semifinite harmonic function on $\zig(t)$ such that $\varphi_{t,w}\geq \psi$. The ideal $\zig(t)\backslash J(t)$ is isomorphic to the $2$-dimensional Pascal graph $\Pas_2$ and the restriction of $\varphi_{t,w}$ to this ideal is an indecomposable harmonic function. Thus, $\varphi_{t,w}$ and $\psi$ are proportional on $\zig(t)\backslash J(t)$ as desired.
\end{example}

\begin{theorem}\label{main theorem}
\phantom{}
\begin{enumerate}[leftmargin=3ex, label=\theenumi)]
    \item For any semifinite zigzag growth model $(t,w)$ the function $\varphi_{t,w}$ is a semifinite indecomposable harmonic function on $\zig(t)$. 
    
    \item Any strictly positive semifinite indecomposable harmonic function on the graph $\zig(t)$ is proportional to $\varphi_{t,w}$ for some semifinite zigzag growth model $(t,w)$.
    
    \item The functions $\varphi_{t,w}$ are distinct for distinct semifinite zigzag growth models $(t,w)$.
\end{enumerate}
\end{theorem}

To prove this theorem we need the following lemma.

\begin{lemma}\label{lemma lemma}
Let $\varphi$ be a strictly positive harmonic function on $\zig(t)$. Then for any $\la\in J(t)$ we have $\varphi\lr{\la}=+\infty$. 
\end{lemma}
\begin{proof}
Recall that 
$$\bw(J(t))=\ml{\ml{\{}}b_0\sqcup \bw(\la^{(1)}) \sqcup b_1\sqcup \ldots \sqcup\bw(\la^{(k)})\sqcup  b_k\  \ml{\ml{\mid}}\ b_i\leq a_i\ \text{and}\ \exists j\colon b_j<a_j;\ \la^{(i)}\in\zig(t_i)\ml{\ml{\}}}.$$
The set on which a harmonic function takes the value $+\infty$ is a coideal. Therefore, without loss of generality, we may assume that
\begin{equation}\label{11122121}
\bw(\la)=a_0\sqcup\bw(\la_1)\sqcup a_1\sqcup\ldots\sqcup\bw(\la_j)\sqcup b_j\sqcup\bw(\la_{j+1})\sqcup a_{j+1}\sqcup\ldots\sqcup\bw(\la_k)\sqcup a_k    
\end{equation}
for some $j$, where $b_j\nearrow a_j$ and each $\bw(\la_i)\in\bw\lr{\zig(v_i)}$ contains as many blocks as possible and all these blocks are large enough. 

We will consider two cases:
\begin{enumerate}[label=\theenumi)]
    \item the block of $a_j$ from which we remove a symbol to obtain $b_j$ is not an outermost block of $a_j$ or it consists of $2$ or more symbols;
    
    \item we remove a symbol from an outermost block of $a_j$ which is of length $1$.
\end{enumerate}

\textit{The first case.} We denote by $\ol{t}$ the template obtained from $t$ by replacing $a_j$ with $b_j$. Then $\la\in\zig(\ol{t})$ and we can write
\begin{equation}\label{1121}
    \biglr{\zig(t)}^{\la}=\biglr{\zig(\ol{t})}^{\la}\sqcup \biglr{\zig(t)}^{\la'},
\end{equation}
where $\la'$ is defined as follows
$$\bw(\la')=a_0\sqcup\bw(\la_1)\sqcup a_1\sqcup\ldots\sqcup\bw(\la_j)\sqcup a_j\sqcup\bw(\la_{j+1})\sqcup a_{j+1}\sqcup\ldots\sqcup\bw(\la_k)\sqcup a_k.$$
Superscripts $\la$ and $\la'$ in \eqref{1121} mean that we take all the zigzags that are greater than or equal to $\la$ and $\la'$ respectively. 

Note that there is an isomorphism of graded graphs 
$$\biglr{\zig(\ol{t})}^{\la}\rightarrow \biglr{\zig(t)}^{\la'},$$
which adds thе removed symbol back $b_j\mapsto a_j$. This map is indeed an isomorphism, since each of these graphs is isomorphic to the Pascal graph of an appropriate dimension. Finally, from Lemma \ref{rem boyers classic} it follows that $\varphi\lr{\la}=+\infty$.

\textit{The second case.} We may assume that $a_j=+-$ or $a_j=-+$. The point is that for all other $a_j$ we can find a vertex which majorizes $\la$ and satisfies the conditions of the first case above. We will restrict ourselves to the case $a_j=-+$. Then we may rewrite \eqref{11122121} as

$$\bw(\la)=a_0\sqcup\bw(\la_1)\sqcup a_1\sqcup\ldots\sqcup\bw(\la_j)\sqcup\bw(\la_{j+1})\sqcup\ldots\sqcup\bw(\la_k)\sqcup a_k$$ 
for some $\la_i\in\zig(t_i)$, because when we delete a symbol from $-+$ inside $\bw(\la)$ the remaining symbol is merged into $\bw(\la_j)$ or $\bw(\la_{j+1})$. Since each $\bw(\la_i)$ contains as many blocks as possible and their lengths are large enough, it follows that $\bw(\la_j)=\al\sqcup \p{n}$ for some "large enough" binary word $\al$ and a natural number $n$. 

Let us define a zigzag $\la'$ as follows
\begin{multline}
\bw(\la')=a_0\sqcup \bw(\la_1) \sqcup a_1\sqcup \ldots \sqcup \bw(\la_{j-1})\sqcup a_{j-1}\\
\sqcup \al\sqcup \p{n-1}\sqcup \m{1}\sqcup\p{1}\sqcup \bw(\la_{j+1})\sqcup a_{j+1}\sqcup\ldots\sqcup\bw(\la_{k-1})\sqcup a_k.    
\end{multline} 

Now we are ready to apply Lemma \ref{gen boyer's lemma}. Below we use the notation from this lemma.  

Let $I\subset \zig(t)$ be the ideal corresponding to the binary words that contain as many blocks as possible, provided that all blocks corresponding to the zigzag flange of $t$ are of maximal lengths. Obviously, this ideal $I$ is isomorphic as a graded graph to the Pascal graph of an appropriate dimension and $\la'\in I$, since $\al$, which appeared in the definition of $\la'$, is "large enough". Then $\dim\lr{\la',\eta}$ on the right hand side of inequality \eqref{ineq gen boyers} from Lemma \ref{gen boyer's lemma} is a multinomial coefficient which arguments are merely differences between lengths of the blocks of $\bw(\la')$ and $\bw(\eta)$ corresponding to infinite clusters of $t$. Then it is easy to check that inequality \eqref{ineq gen boyers} from Lemma \ref{gen boyer's lemma} is fulfilled, because there is a summand on the left hand side of \eqref{ineq gen boyers} which is equal to $\dim\lr{\la',\eta}$. This summand corresponds to $\mu$ obtained from $\eta$ by removing $-$ from $a_j$, appeared in the decomposition of $\bw(\eta)$ provided by Lemma \ref{lemma injection}.
\end{proof}

\begin{proof}[Proof of Theorem \ref{main theorem}]
Let us prove the first two parts of the theorem. We will do it by a single argument. 

Lemma \ref{lemma lemma} implies that any strictly positive harmonic function on $\zig(t)$ is not finite. Then Theorem \ref{theorem ext} provides a bijection between strictly positive indecomposable semifinite harmonic functions on $\zig(t)$ and strictly positive indecomposable finite and semifinite harmonic functions on $\zig(t)\backslash J(t)$. Recall that this bijection is defined by restriction of a function from $\zig(t)$ to $\zig(t)\backslash J(t)$. By Lemma \ref{lemma injection} there is an injective homomorphism of graded graphs 
$$\zig(t)\backslash J(t)\hookrightarrow \zig(t_1)\times\ldots\times\zig(t_k),$$
which image is an ideal. So, applying Theorem \ref{theorem ext} once again we obtain a bijection between strictly positive indecomposable finite and semifinite harmonic functions on $\zig(t)\backslash J(t)$ and 
$$\zig(t_1)\times\ldots\times\zig(t_k).$$ 

Propositions \ref{product graph} and \ref{kernels} yield that this graph admits a finite strictly positive indecomposable harmonic function. Then by Observation \ref{observ fin type2} and Proposition \ref{cunning Boyer} it does not possess any strictly positive semifinite indecomposable harmonic functions. Hence the strictly positive indecomposable finite harmonic functions on this graph are in bijection with strictly positive finite and semifinite indecomposable harmonic functions on $\zig(t)\backslash J(t)$. In particular, this means that the graph $\zig(t)\backslash J(t)$ does not possess strictly positive semifinite indecomposable harmonic functions, because the aforementioned bijection is defined by the restriction of functions from the whole graph to an ideal which is isomorphic to $\zig(t)\backslash J(t)$. 

Next, by Proposition \ref{product graph}, Theorem \ref{theorem fin charr}, and Proposition \ref{kernels} each strictly positive indecomposable finite harmonic function on the graph
$$\zig(t_1)\times\ldots\times\zig(t_k)$$
is of the form 
$$(\la(1),\ldots,\la(k))\mapsto F_{\la(1)}(v_1)\ldots F_{\la(k)}(v_k),$$ 
where $v_i$ is a tuple of oriented consecutive subintervals of $(0,1)$ which orientations are defined by the signs of infinite clusters of $t_i$; the total length of all intervals of all $v_i$'s equals $1$. 

Then each strictly positive finite indecomposble harmonic function on $\zig(t)\backslash J(t)$ is of the form 
$$\la\mapsto  F_{\la^{(1)}}(v_1)\ldots F_{\la^{(k)}}(v_k),$$ 
where $\la\mapsto(\la^{(1)},\ldots,\la^{(k)})$ is the map provided by Lemma \ref{lemma injection}. 

The only thing we are left to do is to indicate how each of these functions should be extended from $\zig(t)\backslash J(t)$ to the whole graph $\zig(t)$. For that purpose we use Lemma \ref{lemma lemma}.

To see that the functions $\varphi_{t,w}$ are distinct for distinct $(t,w)$ we note that the functions 
$$(\la(1),\ldots,\la(k))\mapsto F_{\la(1)}(v_1)\ldots F_{\la(k)}(v_k)$$
are distinct as functions on
$$\zig(t_1)\times\ldots\times\zig(t_k)$$
and by Theorem \ref{theorem ext} the restriction of a function from this graph to an ideal isomorphic to $\zig(t)\backslash J(t)$ is a bijection between finite strictly positive indecomposable harmonic functions.
\end{proof}

\subsection{Semifinite analog of the Vershik-Kerov ring theorem for the zigzag graph}

Recall that by a semifinite analog of the ring theorem we mean Theorem \ref{multiplicativity theorem}. The following proposition describes the finite indecomposable harmonic function that is related to $\varphi_{t,w}$ by this theorem.

\begin{proposition}\label{prop mult for zig}
Let $(t,w)$ be a semifinite zigzag growth model. For any $\mu\in\zig(t)\backslash J(t)$ and $\la\in\zig$ we have 
$$\varphi_{t,w}\lr{F_{\la}F_{\mu}}=\varphi_w(F_{\la})\varphi_{t,w}(F_{\mu}),$$ 
where $\varphi_w$ is the finite harmonic function associated to the finitary oriented paintbox $w$, see Remark \ref{rem456}.
\end{proposition}

Before we start proving this proposition, we need to discuss one preparatory statement, which might be interesting in itself. 

Recall the notation. The zigzag flange of $t$ is denoted by $\fl(t)=(a_0,\ldots,a_k)$ and the splitting of $t$ into sections looks like
$$t=(a_0,\ t_1,\ a_1,\ \ldots\ , a_{k-1},\ t_k,\ a_k).$$

The induced splitting of $w$ reads as $w=v_1\sqcup\ldots\sqcup v_k$.

Let us denote by $v_i(t_i)$ the collection of adjacent intervals on the real line corresponding to the finite template $t_i$ and the tuple of numbers $v_i$. Namely, the length of the $j$-th interval in $v_i(t_i)$ equals $(v_i)_j$ and the orientation of this interval coincides with the sign of the $j$-th infinite cluster of $t_i$. It will be of no importance to us where the leftmost boundary point of $v_i(t_i)$ is placed; only lengths of intervals and their orientations matter.

Let $\eps$ be a real positive number. For any binary word $a$ we denote by $\eps(a)$ a tuple of adjacent oriented intervals each of which corresponds to a block of $a$; orientation of the interval is equal to sign of the block; all intervals are of length $\eps$.

Then we define $w_{\eps}$ as a collection of adjacent intervals from $\eps(a_0),\ldots, \eps(a_k)$ and $v_1(t_1),\ldots, v_k(t_k)$ taken in the order proposed by the splitting of $t$ into sections, that is

$$w_{\eps}=\biglr{\eps(a_0),v_1(t_1),\eps(a_1),v_2(t_2),\ldots,v_k(t_k),\eps(a_k)},$$
see Figure \ref{examp_eps}.

We do not specify where to place the leftmost boundary point of $w_{\eps}$ because it is not important in what follows. Note that the total length of $w_{\eps}$ equals $1+O(\eps)$.

\begin{figure}[h]
    \centering
    \import{Pictures}{Picture_eps-examp.tex}
    \caption{The collection of adjacent intervals $w_{\eps}$ for the semifinite zigzag growth model $(t,w)$ with the semifinite template $\\$ \phantom{paintbox $(t,w)$ with t}$t=\m{1}\pp\mm\p{1}\m{1}\pp\m{2}\pp\m{1}\p{1}\m{2}\pp\mm\p{1}\mm \\$ from Figure \ref{infzig0}. The length and orientation of an interval are indicated above the interval.}
    \label{examp_eps}
\end{figure}

Furthermore, we can define a template $t_{w_{\eps}}$ in the same way as for oriented paintboxes, see the paragraph above Proposition \ref{kernels}. Then $t_{w_{\eps}}$ is a finite template and we may view $w_{\eps}$ as an infinite zigzag $\z(t_{w_{\eps}})$ endowed with a tuple of real positive numbers, see Remark \ref{rm2234432}. Some of these real numbers equal $\eps$ while others come from $w$. The infinite rows and columns to which the number $\eps$ is assigned correspond to the zigzag flange of $t$. So, $\z(t_{w_{\eps}})$ is obtained from $\z(t_w)$ by enlarging the rows and columns corresponding to the binary words from the zigzag flange. Namely, we replace these finite rows and columns with infinite ones, see Figure \ref{fig10}. 

\begin{figure}[h]
    \centering
    \import{Pictures}{Picture_examp10.tex}
    \caption{}
    \label{fig10}
\end{figure}

Let us denote by $\nu_t$ the zigzag obtained from $\z(t)$ by replacing each infinite row or column with a row or column of length $2$. Then $\bw(\nu_t)$ is the binary word obtained from $t$ by replacing each infinite symbol with a single symbol of the same sign.

Obviously, $\nu_t\in\zig(t)\backslash J(t)$ and $\zig(t)^{\nu_t}$, $\biglr{\zig(t_{w_\eps})}^{\nu_t}$ are isomorphic to Pascal graphs of an appropriate dimensions. Recall that the superscript $\nu_t$ means that we take all zigzags greater than or equal to $\nu_t$. 

\begin{proposition}\label{semifin ring th}
There exists a natural number $n$, which depends only on $t$, such that for any $\mu\in\biglr{\zig(t_{w_\eps})}^{\nu_t}$ we have 
\begin{equation}\label{fo}
\varphi_{t,w}\lr{\mu}=\const\cdot \limm_{\eps\to 0}\cfrac{1}{\eps^n}F_{\mu}\lr{w_{\eps}},
\end{equation}
where $\const$ does not depend on $\mu$, but may depend on $(t,w)$; recall that $F_{\la}(w_{\eps})$ is defined by Kerov's construction \eqref{def fin char}, see Section \ref{useful lemma}.
\end{proposition}
\begin{proof}
To obtain the claim we should apply Lemma \ref{lemma for multiplic} to each multiple in Definition \ref{defdefdeff} and to the right hand side of \eqref{fo}.
\end{proof}

\begin{remark}
Equality \eqref{fo} obviously holds for any $\mu\notin\zig(t_{w_\eps})$, because in that case it turns into the trivial identity $0=0$.
\end{remark}

\begin{proof}[Proof of Proposition \ref{prop mult for zig}]

Let us denote by $c_{\la,\mu}^{\nu}$ the structure constants of multiplication in $\qsym$ written in the basis of fundamental quasisymmetric functions, i.e.
$$F_{\la}F_{\mu}=\summ_{\nu}c_{\la,\mu}^{\nu}F_{\nu}.$$

Note that $c_{\la,\mu}^{\nu}\geq 0$ and $c_{\la,\mu}^{\nu}\neq 0$ only if $\nu>\la,\mu$, see \cite[p.35, (3.13)]{quaisymmetric_book}.

Then for any $\la\in\zig$ and $\mu\in\zig(t)^{\nu_t}\subset \biglr{\zig(t_{w_\eps})}^{\nu_t}$ we can write

\begin{multline}
    \varphi_{t,w}\lr{F_{\la}F_{\mu}}=\summ_{\nu}c_{\la,\mu}^{\nu}\varphi_{t,w}\lr{F_{\nu}}=\summ_{\nu}c_{\la,\mu}^{\nu}\const\cdot \limm_{\eps\to 0}\cfrac{1}{\eps^n}F_{\nu}\lr{w_{\eps}}=\\
    \const\cdot \limm_{\eps\to 0}\cfrac{1}{\eps^n}\left[\summ_{\nu}c_{\la,\mu}^{\nu}F_{\nu}\lr{w_{\eps}}\right]=\const\cdot \limm_{\eps\to 0}\left[\cfrac{1}{\eps^n}F_{\la}\lr{w_{\eps}}F_{\mu}\lr{w_{\eps}}\right]=\\
    \const\cdot \left[\, \limm_{\eps\to0}F_{\la}\lr{w_{\eps}}\right]\cdot\left[\,\limm_{\eps\to 0}\cfrac{1}{\eps^n}F_{\mu}\lr{w_{\eps}}\right]=F_{\la}(w)\varphi_{t,w}\lr{F_{\mu}}.
\end{multline}
    
We used formula \eqref{fo} for $\varphi_{t,w}\lr{F_{\nu}}$, because $\nu>\mu$ and then either $\nu\in\biglr{\zig(t_{w_\eps})}^{\nu_t}$ or $\nu\notin\zig(t_{w_\eps})$.

Thus, the finite indecomposable harmonic function $\psi$ from Theorem \ref{multiplicativity theorem} applied to $\varphi_{t,w}$ equals $\varphi_w$.
\end{proof}

    \printbibliography
\end{document}